\documentclass[journal]{IEEEtran}

\usepackage{graphicx}
\usepackage[labelsep=period]{caption}%Fig. 1:, Table 1: -> Fig. 1., Table 1.
\usepackage[caption=false,font=footnotesize]{subfig}
	\usepackage{float}
\usepackage{booktabs}
\usepackage[ruled,vlined]{algorithm2e}
\usepackage{enumerate}
\usepackage{comment}
\usepackage{cite}
\usepackage{color}
\usepackage{xcolor}
\usepackage{hyperref}%把原文中的ref变成跳转链接
	\hypersetup{
		colorlinks,
		linkcolor={green!30!black},
		citecolor={blue!50!black},
		urlcolor={blue!80!black}
	}

\usepackage{amsmath}
\usepackage{amsfonts}
\usepackage{amssymb}
\usepackage{amsthm}
	\newtheorem{lem}{Lemma}
	\newtheorem{thm}{Theorem}
	
	\newtheorem{cor}{Corollary}
	\newtheorem{rmk}{Remark}

	\DeclareMathOperator{\st}{s.t.}

\begin{document}
\bstctlcite{my:BSTcontrol}

\newcommand{\mU}{\text{U}}
\newcommand{\mD}{\text{D}}
\newcommand{\mF}{\text{F}}
\newcommand{\mf}{\text{f}}
\newcommand{\mb}{\text{b}}
\newcommand{\mL}{\text{L}}
\newcommand{\mDG}{\text{DG}}
\newcommand{\mDR}{\text{DR}}
\newcommand{\mE}{\text{E}}

\title{
Closed-Loop Hierarchical Operation for Optimal Unit Commitment and Dispatch in Microgrids: A Hybrid System Approach
%Considering Carbon Emissions Trading
%Closed-Loop Optimal Solution for Unit Commitment in Microgrids Considering Carbon Emissions Trading
%A Hybrid System Approach for Finding the Closed-Loop Optimal Solution of Unit Commitment Problem in Microgrids Considering Carbon Emissions Trading
%    \thanks{%
%    This work was supported in part by the National Natural Science Foundation of China under Grant ***.
%    }
}

\author{
	Fangyuan Li, \IEEEmembership{Student Member, IEEE},
	Jiahu Qin, \IEEEmembership{Senior Member, IEEE},
	%	Wei Xing Zheng, \IEEEmembership{Fellow, IEEE},\\
	and
	Yu Kang, \IEEEmembership{Senior Member, IEEE}

	\thanks{%
		F. Li, J. Qin, and Y. Kang are with the Department of Automation, University of Science and Technology of China,
		Hefei 230027,
		China
		(e-mail:\texttt{fylimas@mail.ustc.edu.cn}; \texttt{jhqin@ustc.edu.cn};\texttt{kangduyu@ustc.edu.cn}).
	}
%	\thanks{%
%		W.~X.~Zheng is with the School of Computing, Engineering and
%		Mathematics, Western Sydney University,
%		Sydney, NSW 2751,
%		Australia
%		(e-mail: \texttt{w.zheng@westernsydney.edu.au}).
%	}
}

\maketitle

\begin{abstract}
To contribute to a smart grid, the unit commitment and dispatch (UCD) is a very important problem to be revisited in the future microgrid environment.
In this paper, the UCD  problem is studied by taking account of both traditional thermal units and distributed generation units (DGs) while considering demand response (DR) and carbon emissions trading (CET).
Then, the UCD problem is reformulated as the optimal switching problem of a hybrid system.
The reformulation allows convenient handling of time-dependent start-up cost for dynamic programming (DP).
Afterwards, a closed-loop hierarchical operation (CLHO) algorithm is developed which is capable of producing the optimal schedule in the closed-loop behavior.
Theoretical analysis is also presented to show the equivalence of the reformulation and the effectiveness of the CLHO algorithm.
Finally, the closed-loop behavior and the effectiveness of the algorithm are further verified by simulation studies.
It is also shown that low emission quotas and high emission trading price help reduce the total carbon emissions.
\end{abstract}


\section{Introduction\label{sec:Intro}}
The microgrid is a promising concept to help integrating various distributed generation units (DGs) in smart grid. % \cite{microgrid}. 
To cope with the new factors such as DGs and demand response (DR) in microgrid, the unit commitment and dispatch (UCD) is an important problem that should be revisited. 
The objective of UCD is to schedule the generation units such that the total generation cost over a time horizon is minimized while satisfying power balance constraints and various operational constraints \cite{Wood:2013wo}.

Traditional methods for solving UCD problem include priority-list (PL) \cite{UCviaPL}, Lagrangian relaxation (LR) \cite{UCviaLR}, and various heuristic search algorithms \cite{UCviaPSO}.
In PL method, the shutdown rule or priority-list scheme acts a key role in deciding which unit to commit or drop in order to satisfy the increasing or dropping demand. The PL methods are usually simple and fast but result in relatively high total generation costs.
The LR method is based on using Lagrange multipliers to incorporate the constraints onto the total cost, which results in the so-called Lagrangian. The UCD problem is converted to find the primal and dual solution of the Lagrangian. It is relatively easy to add constraints in the Lagrangian for LR method while the main disadvantage is the inherent suboptimality \cite{Wood:2013wo}.
The dynamic programming (DP) method is also widely employed in solving UCD problem. The distinct advantage of DP method is the ability to find the optimal schedule for the UCD problem. However, it is generally difficult for the DP method to deal with the time-dependent start-up costs and constraints \cite{Padhy:2004eq}.
The heuristic search methods including particle swarm algorithm, evolution programming, and genetic algorithm have been continuously studied to solve UCD problem. The main advantage of these methods is the capability of solving large scale problems under various constraints. However, there is no general rule of tuning the parameters of the algorithms \cite{UCviaPSO}.
It is worth noting that though the UCD problem has been extensively studied in various methods, the problem remains challenging due to the nonlinearity and the presence of both integer and continuous variables.

The recent advances in reinforcement learning (RL) \cite{Sutton:1998wc,Littman:2015hm} and DP \cite{Bertsekas:2013wv,Kiumarsi:2018kq} areas are providing new promising methods for solving some optimization problems in power systems.
A series of successful applications have been made to solve various economic power dispatch problems (ED) which contain no integer variables, see \cite{Wei:2017hp,Yang:2018km,Zhu:2018kj,Shuai:2018et} for examples.
However, %comparing to various modified heuristic search algorithms\cite{UCviaPSO1,UCviaPSO2,Zheng:2015fq}, 
relatively little works exploit these advances to solve the more challenging UCD problem where integer variables are also involved.
In \cite{Xie:2009ju}, model predictive control is applied to solve the UCD problem with intermittent resources while the start-up and shutdown costs are not considered.
Reference \cite{Zhang:2011dg} employs state-space approximate dynamic programming (ADP) for solving the stochastic UCD under the framework of Markov decision process. However, only the commitment status is treated as decision variables in the UCD problem without taking account of power generations.
Recently, reference \cite{iconip} proposes a consensus based distributed reinforcement learning algorithm to solve the non-convex ED problem considering valve-point effects. The algorithm provided in \cite{iconip} can also be extended to solve the UCD problem.
However, a drawback of \cite{iconip} is the requirement of the discrete generation output which may introduce discretization errors.
Nevertheless, these works have not taken account of the effect of the DR and carbon emissions trading (CET).

In this paper, the deterministic UCD problem is investigated in the microgrid environment consisting of thermal units and DGs while considering DR and CET. 
%The UCD problem is reformulated as the optimal switching problem of a hybrid system.
%Then, a closed-loop hierarchical operation (CLHO) algorithm is developed using neural networks to approximate the optimal cost-to-go of the hybrid system and using convex optimization algorithm to solve the lower level problem.
The contributions are listed as follows.
1) Besides the traditional thermal units, the DGs, DR, and CET are also considered. With the increasing penetration of renewable generation %, research interest of DR, 
and pressure of emissions reduction, the UCD problem considered in this paper makes a more practical scenario than existing works such as \cite{Xie:2009ju,Zhang:2011dg,iconip}. 
2) The UCD problem is firstly reformulated as the optimal switching problem of a hybrid system. 
In the reformulation, the time-dependent start-up cost is converted to the switching cost of the hybrid system, enabling the dealt of this cost which is generally seen as difficult for DP \cite{Padhy:2004eq}.
The reformulation also allows exploiting advances in related areas % such as RL and DP 
\cite{Sutton:1998wc,Littman:2015hm,Bertsekas:2013wv,Kiumarsi:2018kq} to solve the UCD problem. 
3) A closed-loop hierarchical operation (CLHO) algorithm is developed to find the optimal schedule of the hybrid system.
The proposed CLHO algorithm is capable of producing the optimal schedule in a closed-loop style, where the trained neural network can deal with different initial conditions and disturbances without retraining. In contrast, most existing methods are usually open-loop, thus cannot guarantee optimality under external disturbances without reoptimization. Theoretical analysis 
	%on both the equivalence of the reformulated problem and the effectiveness of the proposed algorithm 
	are also provided as opposed to most existing works such as \cite{Xie:2009ju,Zhang:2011dg,Rong:2017ip}.
4) Simulation results verify that the CLHO algorithm produces the optimal schedule in closed-loop under different initial conditions and external disturbance, which potentially helps the real-time optimal operation of the microgrids. It is also shown in the simulation that low emission quotas and high emission trading price reduce the total amount of carbon emissions.

\begin{comment}
%贡献和创新点
1. 考虑碳排放，是的微电网模型更加完善，相比于已有文献更具有挑战性；
2. 将问题形式化为混杂系统模型，并给出等价性条件，为使用新方法求解提供条件；
3. 针对这一问题，我们给出了一种新的闭环求解算法，并提供了严格的分析；
4. 设置了考虑碳排放的微电网仿真系统。验证了提出方法的优越性。
\end{comment}

The paper is organized as follows.
Section \ref{sec:Pre} provides some preliminaries. % on microgrid environment and CET.
The UCD problem is reformulated as an optimal switching problem of a hybrid system in Section \ref{sec:Pro}. 
Section \ref{sec:Alg} develops the CLHO algorithm along with theoretical analysis. % to guarantee the effectiveness.
In Section \ref{sec:Sim}, simulations are conducted to study the behavior and effectiveness of the proposed algorithm.
Finally, we draw conclusions in Section \ref{sec:Con}.

\section{Preliminary\label{sec:Pre}}
In this section, some preliminaries on CET and microgrid environment are presented.

\subsection{Carbon Emissions Trading}
CET is a cap-and-trade system where the government sets a cap of carbon emissions and issues quotas for each emitter.
The quotas can be seen as allowances for emitters to emit certain amount of greenhouse gas.
The emitters can trade quotas in the CET market as needed, i.e., sell surplus or purchase deficit of emission allowances.
In the UCD problem in this paper, the emitters are thermal units.
As argued in paper \cite{Zhang:2016cl}, we assume that the trade of generation units in UCD does not affect the trading price.
The emission cost of thermal units is determined by
\begin{equation}
C_{E} = \sum_{n=1}^{N} \left(\sum_{t=1}^{T} E_n (P_n[t]) I_n[t] - Q_{n} \right)\cdot p_e,
\label{eq:ce}
\end{equation}
\begin{equation}
E_n (P_n[t]) = \alpha_n (P_n[t])^2 + \beta_n P_n[t] + \gamma_n,
\end{equation}
where $E_n (P_n[t])$ is the emission of thermal unit $n$ at period $t$, $Q_n$ is the emission quota of unit $n$, $p_e$ is the CET price in the market, $\alpha_n$, $\beta_n$, and $\gamma_n$ are the emission coefficients.

\subsection{Microgrid Environment}
\subsubsection{Thermal Units}
The fuel cost of thermal unit $n$ is usually modeled as a quadratic function \cite{advances}:
\begin{equation}
C_{\mF,n}(P_n[t])=a_n (P_n[t])^2 + b_n P_n[t] + c_n,
\end{equation}
where $P_n[t]$ is the active power output of unit $n$ at period $t$, $a_n$, $b_n$, and $c_n$ are the cost coefficients of unit $n$.

To bring the thermal unit online, a certain amount of energy must be spent, incurring a start-up cost.
Assume that the off-line thermal unit $n$ is in banking mode, then the start-up cost is usually modeled as \cite{Wood:2013wo}
\begin{equation}
C_{\mU,n}(\tau_n)=C_{\mb,n} \cdot \tau_n + C_{\mf,n},
\label{eq:cu}
\end{equation}
where $C_{\mb,n}$ is the cost for unit $n$ to maintain the operating temperature, $\tau_n$ is the period that unit $n$ was in banking mode, and $C_{\mf,n}$ represents the fixed cost.
To decommit a thermal unit, a fixed shutdown cost $C_{\mD,n}$ is usually incurred.

Define $I_n[t]$ as the commitment decision variable of unit $n$ at time $t$, where $I_n[t]=1$ if unit $n$ is committed/online at $t$ and $I_n[t]=0$ if otherwise.
The power output of unit $n$ usually subjects to the following generation capacity constraints:
\begin{equation}
I_n[t] P^m_n \leq P_n[t] \leq I_n[t] P^M_n,
\end{equation}
where $P^m_n$ ($P^M_n$) is the minimum (maximum) power output of unit $n$.

The power output variation of thermal units usually subject to the following ramp limit constraint
\begin{equation}
-R_{n}^{\mD} \leq (P_{n}[t]-P_{n}[t-1]) I_{n}[t] I_{n}[t-1] \leq R_{n}^{\mU},
\end{equation}
where $R_{n}^{\mD}$ and $R_{n}^{\mU}$ are the ramp down and ramp up rate limit.%, respectively.

\subsubsection{Distributed Generation}
Suppose there is an aggregator that manages the DGs such as wind turbines and solar panels, then the aggregator can be treated as a virtual generation unit with cost function \cite{Zhang:2016cl}
\begin{equation}
C_{\mDG}(P_{\mDG}[t]) = a_{\mDG} (P_{\mDG}[t])^2 + b_{\mDG}P_{\mDG}[t] + c_{\mDG},
\end{equation}
where $a_{\mDG}$, $b_{\mDG}$, and $c_{\mDG}$ are the cost coefficients, and $P_{\mDG}[t]$ is the aggregated generation.

Since that DGs subject to natural conditions and tend to be intermittent, upper limits on available DGs and their penetration rate %at period $t$ 
are usually imposed as following \cite{Zhang:2016cl}, to ensure a reliable operation
\begin{equation}
0\leq P_{\mDG}[t]\leq P^{M}_{\mDG}[t],
\end{equation}
and
\begin{equation}
\eta_{t} = \frac{P_{\mDG,t}}{\sum_{n=1}^N P_n[t] + P_{\mDG,t}}\leq \eta_{\max}.
\end{equation}
\subsubsection{Demand Response}
DR captures the relation between the load elasticity and the electricity price.
The demand curtailment due to DR can be treated as a virtual generation unit \cite{Aghaei2013}.
The cost function of this virtual generation unit can be modeled as follows \cite{Aghaei2013}
\begin{equation}
C_{\mDR}(P_{\mDR}[t]) = a_{\mDR}(P_{\mDR}[t])^2 + b_{\mDR} P_{\mDR}[t] + c_{\mDR},
\end{equation}
where $a_{\mDR}$, $b_{\mDR}$, and $c_{\mDR}$ are the cost coefficients, and $P_{\mDR}[t]$ is the curtailed demand.

There is an upper limit on the demand curtailment at each period $t$ as follows
\begin{equation}
0\leq P_{\mDR}[t] \leq P^{M}_{\mDR}[t].
\end{equation}

\subsubsection{Operational Constraints}
The operation of microgrid subject to the power balance constraints \cite{yanni}
\begin{equation}
\sum_{n=1}^N P_n[t]I_n[t] + P_{\mDG}[t] + P_{\mDR}[t] = P_{\mD}[t], \quad t=1,...,T.
\label{eq:balance}
\end{equation}
where $P_{\mD}[t]$ is the power demand, $N$ is the total number of thermal units, and $T$ is the total number of time periods.
%The power losses $P_{\mL}[t]$ is usually modeled by Kron's loss formula $\sum_{i=1}^N\sum_{j=1}^NP_i[t]B_{ij}P_j[t]+\sum_{i=1}^N B_{0i}P_i[t] + B_{00}$, where $B_{ij}$, $B_{0i}$, and $B_{00}$ are loss coefficients.
Since an isolated microgrid environment is considered in this paper, there is no power exchange with the main grid in \eqref{eq:balance}.

To cope with the fluctuation of renewable energies and demand, the following spinning reserves are usually considered
\begin{equation}
\sum_{n=1}^N P_n^m I_n[t] + P_{\mDG}[t] + P_{\mDR}[t] - P_{\mD}[t]\leq -R^m[t],
\end{equation}
\begin{equation}
\sum_{n=1}^N P_n^M I_n[t] + P_{\mDG}[t] + P_{\mDR}[t] - P_{\mD}[t]\geq R^M[t],
\end{equation}
where $R^m[t]>0$ and $R^M[t]>0$ are the lower and upper spinning reserves.
\subsubsection{Power Flow Constraints}
Note that the primary focus of this paper is solving the UCD problem, although the optimal power flow problem is also related.
Hence, the transmission losses and power flow constraints are not incorporated  explicitly in the following UCD model.
See \cite{Molzahn:2017ba} for various power flow representations.
However, it is worth noting that both the power flow constraints (linear approximations or convex approximations) and the transmission losses can be added.
The proposed algorithm is still valid with appropriate decision variables since the lower level problem remains convex. % considering transmission losses and power flow constraints with linear or convex approximations.

\section{Problem Formulation\label{sec:Pro}}
In this section, the UCD problem is presented in a commonly used form as in literature \cite{Padhy:2004eq,Zhang:2016cl}.
Then, the problem is decomposed and reformulated into an optimal switching problem of a hybrid system.
%Theoretical results are also provided to ensure the equivalence of the reformulated problem.
\subsection{Unit Commitment and Dispatch Model}
The UCD model considering DR, DG, and CET can be formulated as follows \cite{Wood:2013wo,Padhy:2004eq}:
\begin{align}
\min\ &\sum_{t=1}^T \Big(\sum_{n=1}^N ( C_{\mF,n}(P_n[t])I_n[t] + C_{\mU,n}[t] + C_{\mD,n}[t] )\notag\\
&\qquad\qquad+ C_{\mDG}(P_{\mDG}[t]) + C_{\mDR}(P_{\mDR}[t]) \Big)+ C_{\mE}\label{eq:uc}\\
\st\ &\sum_{n=1}^N P_n[t]I_n[t] + P_{\mDG}[t] + P_{\mDR}[t] = P_{\mD}[t],
\label{eq:pbc}\\
&\sum_{n=1}^N P_n^m I_n[t] + P_{\mDG}[t] + P_{\mDR}[t]\notag\\
&\qquad\leq P_{\mD}[t]- R^m[t],
\label{eq:sr1}\\
&\sum_{n=1}^N P_n^M I_n[t] + P_{\mDG}[t] + P_{\mDR}[t]\notag\\
&\qquad\geq P_{\mD}[t]+ R^M[t],
\label{eq:sr2}\\
&I_n[t] P^m_n \leq P_n[t] \leq I_n[t] P^M_n,
\label{eq:gcc}\\
&-R_{n}^{\mD} \leq (P_{n}[t]-P_{n}[t-1]) I_{n}[t] I_{n}[t-1] \leq R_{n}^{\mU},
\label{eq:rl}\\
&0\leq P_{\mDG}[t]\leq P^{M}_{\mDG}[t],
\label{eq:dgc}\\
&\eta_{t} = \frac{P_{\mDG}[t]}{\sum_{n=1}^N P_n[t] + P_{\mDG}[t]}\leq \eta_{\max}
\label{eq:eta}\\
&0\leq P_{\mDR}[t] \leq P^{M}_{\mDR}[t],
\label{eq:drc}
\end{align}
where the constraints \eqref{eq:pbc}--\eqref{eq:drc} are imposed for all $n=1,...,N$ and $t=1,...,T$. 
%Similar to \cite{Padhy:2004eq} where $t$ is used in subscript, one should keep in mind that $C_{\mU,n}[t]$ represents the possible start-up cost induced at $t$.
$C_{\mU,n}[t]$ is computed by \eqref{eq:cu} if $(1-I_n[t-1])I_n[t]=1$, otherwise, $C_{\mU,n}[t]=0$.
$C_{\mD,n}[t]$ equals $C_{\mD,n}$ if $I_n[t-1](1-I_n[t])=1$, otherwise, $C_{\mD,n}[t]=0$.
It is assumed that all dispatchable thermal units and nondispatchable DGs belong to a utility microgrid.
\begin{comment}
\begin{rmk}
Note that there is a common way of incorporating the start-up and shutdown cost in the objective, namely, in form of $\sum_{t=1}^T C_{\mU,n}(1-I_n[t-1])I_n[t]+C_{\mD,n}I_n[t-1](1-I_n[t])$, see reference \cite{UCviaPSO} and \cite{Zhang:2016cl} respectively for examples.
However, this notation leaves out considering the cost of the units decommitted at the tail of the $T$ horizon when the units are in banking mode.
For example, if a thermal unit is off-line at $[T-1,T]$, start-up cost still rises if in banking mode while the common notation fails to incorporate this cost.
Therefore, by referring to $\sum_{t=1}^T C_{\mU,n}[t] + C_{\mD,n}[t]$ in \eqref{eq:uc}, one should keep in mind that all the possible incurred start-up and shutdown cost at each time $t$ are incorporated.
It should not be confused with $C_{\mU,n}(t)$ where $t$ is the independent variable of $C_{\mU,n}$.
\end{rmk}
\end{comment}

\subsection{Model Decomposition}
Substituting \eqref{eq:ce} into objective \eqref{eq:uc}, one obtains %the following
\begin{multline}
\sum_{t=1}^T \Big( \sum_{n=1}^N \big( I_n[t] (C_{\mF,n}(P_n[t]) + E_n(P_n[t])) -Q_n p_e\big) + C_{\mDG}(P_{\mDG}[t])\\
+ C_{\mDR}(P_{\mDR}[t]) \Big)
+ \sum_{t=1}^T\sum_{n=1}^N(C_{\mU,n}+C_{\mD,n}).
\label{eq:md}
\end{multline}
The first item in \eqref{eq:md} incorporates all the fuel costs, emissions costs, and the costs of DGs and DR of each time period $t$,
while the second item summarizes all the costs incurred by the commitment or decommitment of thermal units along the time horizon $T$.
With this observation, we denote the first and the second item compactly as
\begin{multline}
\sum_{t=1}^T Q(P[t],I[t])
:=\sum_{t=1}^T \Big( \sum_{n=1}^N \big(I_n[t](C_{\mF,n}(P_n[t]) + E_n(P_n[t]))\\
 -Q_n p_e\big) + C_{\mDG}(P_{\mDG}[t]) + C_{\mDR}(P_{\mDR}[t]) \Big)
\label{eq:q}
\end{multline}
and
\begin{equation}
\sum_{t=1}^T \kappa(I[t-1],I[t]) :=\sum_{t=1}^T\sum_{n=1}^N(C_{\mU,n}+C_{\mD,n})
\label{eq:k}
\end{equation}
where $P[t]$ and $I[t]$ are vector $[P_1[t], ..., P_N[t], P_{\mDG}[t], P_{\mDR}[t]]^T$ and $[I_1[t],...,I_N[t]]^T$, respectively.
In \eqref{eq:k}, $I[0]$ represents the initial commitment status of the thermal units.
Then, the objective \eqref{eq:uc} can be denoted as
\begin{equation}
\sum_{t=1}^T Q(P[t],I[t])+\kappa(I[t-1],I[t]).
\label{eq:qk}
\end{equation}
Under certain conditions, the original UCD problem can be seen as a hierarchical optimization problem, where the upper level objective is to minimize \eqref{eq:qk} along %the total time 
$T$ time horizon and the lower level objective is to minimize \eqref{eq:q} subjecting to constraints \eqref{eq:pbc}--\eqref{eq:drc} at each time $t$.
\subsection{Hybrid System Model}
%According to \eqref{eq:cu} and \eqref{eq:k}, 
It can be verified that the start-up and shutdown cost $\kappa(I[t-1],I[t])$ incurred at each time period $t$ 
can be written as
\begin{equation}
\kappa(I[t-1],I[t]) = \sum_{n=1}^N \kappa_n(I[t-1],I[t])
\label{eq:kappa}
\end{equation}
where $\kappa_n(I[t-1],I[t])$ is calculated by
\begin{multline}
\kappa_n(I[t-1],I[t])=
C_{\mb,n}+(C_{\mf,n}-C_{\mb,n}+C_{\mD,n})I_n[t-1]\\
- (C_{\mf,n}+C_{\mD,n})I_n[t-1]I_n[t].
\label{eq:kn}
\end{multline}
The introduction of \eqref{eq:kappa} and \eqref{eq:kn} eliminates the need of an implicit variable $\tau_n$ to keep track of the banking time as in literature \cite{Wood:2013wo,Padhy:2004eq}. 
%In addition, the total start-up and shutdown cost \eqref{eq:k} which requires commitment decisions of multiple $\tau_n$ length can be computed cumulatively with only consequent decisions $I_n[t-1]$ and $I_n[t]$.
%An additional benefit is that the incurred start-up and shutdown cost are incorporated with mathematically rigorous notation.
Since the cost is incurred by switching from $I[t-1]$ to $I[t]$, we call $\kappa(I[t-1],I[t])$ switching cost. %in the following of this paper.

Note that there are two exogenous variables $P^M_{\mDG}[t]$ and $P_{\mD}[t]$ in the UCD problem. %that determine the minimum total cost and optimal solution of the UCD problem.
In the day-ahead UC, the short-time renewable power and load forecasting are usually assumed available, which can be denoted as
\begin{equation}
P^M_{\mDG}[t] = F_1 (\bar{P}^M_{\mDG}[t-1];\theta_1),
\label{eq:rpf}
\end{equation}
and
\begin{equation}
P_{\mD}[t] = F_2 (\bar{P}_{\mD}[t-1];\theta_2),
\label{eq:lf}
\end{equation}
where $F_{(\cdot)}$ are the forecasting models, $\bar{P}^M_{\mDG}[t]$ and $\bar{P}_{\mD}[t]$ are the regression vectors, and $\theta_{(\cdot)}$ are model parameters.
Interested readers are directed to \cite{Espinoza:2007bf} and references therein for details on the establishment of these models.
Then, the lower level optimization problem can be rewritten as
\begin{equation}
\begin{cases}
&P[t]=\arg\min Q(P[t],I[t])\\
&\st\ \eqref{eq:pbc},\eqref{eq:sr1},\eqref{eq:sr2},\eqref{eq:gcc},\eqref{eq:rl},\eqref{eq:dgc},\eqref{eq:eta},\eqref{eq:drc}
\end{cases}.
\label{eq:lop}
\end{equation}
Denote \eqref{eq:lop} as
\begin{equation}
P[t] = f_I (P[t-1],P^M_{\mDG}[t],P_{\mD}[t]),
\label{eq:lop2}
\end{equation}
where subscript $I\in\{0,1\}^N$ is the index of subsystems.
Note that for some $I\in\{0,1\}^N$ with $P_{\mDG}^M[t]$, $P_{\mD}[t]$, $R^m[t]$, $R^M[t]$, and $P^M_{\mDR}[t]$, the lower optimization problem may admit no solutions.
To eliminate these exceptions, we use $\mathcal{I}(P_{\mDG}^M[t],P_{\mD}[t],R^m[t],R^M[t],P^M_{\mDR}[t])$, in shorthand as $\mathcal{I}_t$ for given variable values, to denote the set of $I$ when the lower optimization problem admits a solution $P[t]\in\mathcal{P}_t$. $\mathcal{P}_t$ denotes the set of feasible solutions.

Finally, %considering \eqref{eq:qk} and \eqref{eq:rpf}\eqref{eq:lf}\eqref{eq:lop2},
the original UCD model is reformulated as minimizing the objective 
\begin{equation}
\sum_{t=1}^T Q(P[t],I[t])+\kappa(I[t-1],I[t]).
\label{eq:cost}
\end{equation}
while subjecting to a hybrid system as follows
\begin{equation}
\left\{
\begin{aligned}
&P^M_{\mDG}[t] = F_1 (\bar{P}^M_{\mDG}[t-1];\theta_1),\\
&P_{\mD}[t] = F_2 (\bar{P}_{\mD}[t-1];\theta_2),\\
&P[t] = f_I (P[t-1],P^M_{\mDG}[t],P_{\mD}[t]),
\end{aligned}
\right.
\label{eq:hybrid}
\end{equation}
Note that the hybrid system \eqref{eq:hybrid} has $|\mathcal{I}|$ subsystems where $|\cdot|$ denotes cardinality.
The subsystems are indexed by $I\in \mathcal{I}$.
At each time $t$, only one subsystem $I[t]\in \mathcal{I}$ is running.
The objective is to choose a sequence of $I[t]$ such that \eqref{eq:cost} is minimized.
Interested readers are directed to \cite{Tomlin:2000us,Goebel:2011tx} %\cite{Tomlin:2000us,Goebel:2009gr,Goebel:2011tx}
for more details on hybrid systems.

Since all the start-up and shutdown cost are incorporated in \eqref{eq:cost} without involving any implicit variables, easy handling of these costs is allowed which is generally seen as difficult for DP. 
In the next section, DP and approximation approach will be employed to solve the reformulated problem. 
To this end, we show the equivalence between the reformulated problem and the original UCD model. 
%minimizing the cost \eqref{eq:cost} under the hybrid system \eqref{eq:hybrid} and minimizing the objective \eqref{eq:uc} under constraints \eqref{eq:pbc}--\eqref{eq:drc}.
\begin{comment}
There may exist multiple optimal solutions to an ordinary optimization problem, however it usually requires the system dynamic to be single-valued, otherwise, the system trajectory may not be unique.
It is shown in the following result that the formulated hybrid system \eqref{eq:hybrid} is single-valued and admits a unique trajectory.
\end{comment}
Before that, we have to show that the formulated hybrid system \eqref{eq:hybrid} is single-valued and admits a unique trajectory.
\begin{thm}
Given a (feasible) switching schedule $(I[t])_{t=1}^T$, the hybrid system \eqref{eq:hybrid} admits a unique trajectory $(P[t])_{t=1}^T$.
\label{thm:prop}
\end{thm}
\begin{IEEEproof}
The existence of a switching schedule $(I[t])_{t=1}^T$ implies that $\mathcal{I}_t\neq\varnothing$ for any $t=1,...,T$, which means that the lower level optimization problem \eqref{eq:lop} admits a solution at all $t=1,...,T$.
Thus, the hybrid system admits some trajectory $(P[t])_{t=1}^T$.
It remains to show that the trajectory is also unique.

Denote $P[0]$ and $I[0]$ the initial condition.
Then, $P[1]$ is the optimal solution of problem \eqref{eq:lop} with $t=1$ for $I[1]$ and initial condition $P[0]$ and $I[0]$.
The objective $Q(P[1],I[1])$ is strictly convex since for all $n=1,...,N$, $C_{\mF,n}(P_n[1])$, $E_n(P_n[1])$, $C_{\mDG}(P_{\mDG}[1])$, and $C_{\mDR}(P_{\mDR}[1])$ are strictly convex w.r.t. $P[1]$, and $Q_n p_e$ is constant.
The feasible regions of all the constraints \eqref{eq:pbc}--\eqref{eq:drc} are convex, leading to a convex intersection region.
Then, the uniqueness of $P[1]$ follows. % convexity of problem \eqref{eq:lop} with $t=1$ for $I[1]$.
Similar argument can be applied to any $t=1,...,T$.
This completes the proof.
\end{IEEEproof}

Once the initial $P[0]$ and $I[0]$ are given, and a switching schedule $\{I[t]\}_{t=1}^T$ is selected, the hybrid system can operate from the initial time $t=1$ to the final time $t=T$, incurring running cost $Q(P[t],I[t])$ and switching cost $\kappa(I[t-1],I[t])$.
The objective of the problem is to find an optimal switching schedule that minimizes the total cost \eqref{eq:cost}.

Then we show that the hybrid system \eqref{eq:hybrid} minimizing cost \eqref{eq:cost} is equivalent to the original UCD with objective \eqref{eq:uc} and constraints \eqref{eq:pbc}--\eqref{eq:drc} under certain conditions.
\begin{thm}
Assume that the original UCD problem with objective \eqref{eq:uc} and constraints \eqref{eq:pbc}--\eqref{eq:drc} admits feasible solutions.
Then, the hybrid system \eqref{eq:hybrid} w.r.t. cost function \eqref{eq:cost} admits an optimal switching schedule and a unique corresponding trajectory, denoted by $(I^*[t])_{t=1}^T$ and $(P^*[t])_{t=1}^T$.
Moreover, if the ramp limit constraint \eqref{eq:rl} is relaxed, the optimal switching schedule $(I^*[t])_{t=1}^T$ along with the unique corresponding trajectory $(P^*[t])_{t=1}^T$ minimizes the original UCD problem.
\label{thm:eqv}
\end{thm}
\begin{IEEEproof}
Firstly, we show that the hybrid system admits an optimal switching schedule.
It is sufficient to show the system admits feasible switching schedules.
Then the finiteness leads to the existence of an optimal switching schedule.

Denote the initial conditions by $P[0]$ and $I[0]$, and denote $(P[t])_{t=1}^T$ and $(I[t])_{t=1}^T$ an feasible solution of the original UCD problem.
Then, $P[t]$ and $I[t]$ are also a feasible solution of the lower level optimization problem \eqref{eq:lop} with initial conditions $P[t-1]$ and $I[t-1]$ for any $t=1,...,T$.
By definition of hybrid system \eqref{eq:hybrid}, $(P[t])_{t=1}^T$ is also a trajectory under the switching schedule $(I[t])_{t=1}^T$.
Namely, the hybrid system admits switching schedule $(I[t])_{t=1}^T$.
Then, it follows that the hybrid system admits an optimal switching schedule.

Denote $(I^*[t])_{t=1}^T$ as the optimal switching schedule.
Then, the Theorem \ref{thm:prop} indicates that the optimal switching schedule admits a unique trajectory, denoted by $(P^*[t])_{t=1}^T$.
The above argument still holds if ramp limit constraint \eqref{eq:rl} is relaxed.

Secondly, we show that if ramp limit constraint is relaxed, the optimal schedule along with the unique corresponding trajectory  minimizes the original UCD problem by contradiction.

Assume that the schedule $(I^*[t])_{t=1}^T$ along with trajectory $(P^*[t])_{t=1}^T$ does not minimize the original UCD problem.
Then there exist some $(P'[t])_{t=1}^T$ and $(I'[t])_{t=1}^T$ such that
\begin{multline}
\sum_{t=1}^T Q(P'[t],I'[t]) + \kappa(I'[t-1],I'[t]) \\
< \sum_{t=1}^T Q(P^*[t],I^*[t]) + \kappa(I^*[t-1],I^*[t]).
\label{eq:leq}
\end{multline}
Note that $(P'[t])_{t=1}^T$ does not necessarily minimize the original problem nor the reformulated problem under the schedule $(I'[t])_{t=1}^T$.
If it holds that $(I^*[t])_{t=1}^T = (I'[t])_{t=1}^T$, then \eqref{eq:leq} leads to
$(P^*[t])_{t=1}^T = (P'[t])_{t=1}^T$ by contradiction.
Notice that with \eqref{eq:leq}, $(I^*[t])_{t=1}^T = (I'[t])_{t=1}^T$ leads to
\begin{equation}
\sum_{t=1}^T Q(P'[t],I'[t]) < \sum_{t=1}^T Q(P^*[t],I^*[t]).
\end{equation}
Then, there must exists some $t=\tau\in\{1,...,T\}$ such that
\begin{equation}
Q(P'[t],I'[t]) = Q(P'[t],I^*[t]) < Q(P^*[t],I^*[t]).
\end{equation}
However, this contradicts the fact that $Q(P^*[t],I^*)\leq Q(P'[t],I^*[t])$ for any $P'[t]$ according to \eqref{eq:lop} and \eqref{eq:hybrid} in which the ramp limit constraints have been relaxed.
Thus, it holds that $Q(P'[t],I^*[t]) = Q(P^*[t],I^*[t])$ for all $t\in\{1,...,T\}$, leading to $(P^*[t])_{t=1}^T = (P'[t])_{t=1}^T$ with Theorem \ref{thm:prop}.

Denote $(P''[t])_{t=1}^T$ as the system trajectory operated under switching schedule $(I'[t])_{t=1}^T$.
If otherwise $(I^*[t])_{t=1}^T \neq (I'[t])_{t=1}^T$, then one obtains
\begin{equation}
\sum_{t=1}^T Q(P''[t],I'[t]) \leq \sum_{t=1}^T Q(P'[t],I'[t])
\end{equation}
by observing \eqref{eq:lop} and \eqref{eq:hybrid} where the ramp limit constraints have been relaxed.
Thus, it holds that
\begin{multline}
\sum_{t=1}^T Q(P''[t],I'[t]) + \kappa(I'[t-1],I'[t]) \\
\leq \sum_{t=1}^T Q(P'[t],I'[t]) + \kappa(I'[t-1],I'[t])\\
< \sum_{t=1}^T Q(P^*[t],I^*[t]) + \kappa(I^*[t-1],I^*[t]).
\end{multline}
That is, schedule $(I'[t])_{t=1}^T$ with trajectory $(P''[t])_{t=1}^T$ results in lower total cost than schedule $(I^*[t])_{t=1}^T$ with trajectory $(P^*[t])_{t=1}^T$.
Since $(I^*[t])_{t=1}^T$ has been assumed as an optimal switching schedule, contradiction is raised.
\end{IEEEproof}

\begin{comment}
A byproduct of Theorem \ref{thm:eqv} is %that the original UCD problem admits optimal solutions if ramp limit constraint is relaxed, 
summarized in the following corollary.
\begin{cor}
Assume that the original UCD problem with objective \eqref{eq:uc} and constraints \eqref{eq:pbc}--\eqref{eq:gcc} and \eqref{eq:dgc}--\eqref{eq:drc} admits feasible solutions.
Then, there exists optimal $(P^*[t])_{t=1}^T$ and $(I^*[t])_{t=1}^T$ such that the objective \eqref{eq:uc} is minimized.
\end{cor}
\begin{IEEEproof}
The existence of optimal switching schedule $(I^*[t])_{t=1}^T$ along with the unique corresponding trajectory $(P^*[t])_{t=1}^T$ , as argued in the proof of Theorem \ref{thm:eqv}, immediately proves this result.
Alternatively, note that the continuity of objective functions on the compact feasible set for any feasible schedule also lead to this result. 
\end{IEEEproof}
\end{comment}
Considering the above equivalence, efforts are made to find the optimal switching schedule of hybrid system \eqref{eq:hybrid} w.r.t. the cost function \eqref{eq:cost} in the next section.

\section{Proposed Algorithm\label{sec:Alg}}
In this section, the optimal cost-to-go of the hybrid system is presented first, then the CLHO algorithm is developed where neural networks are employed to approximate the optimal cost-to-go.
%In order to show the validity and effectiveness for solving the optimal switching problem, theoretical analysis is also developed.
\subsection{Theory}
Denote $J_t^{I[t:T]}(P[t],I[t-1])$ as the cost-to-go associated with $P[t]$ and $I[t-1]$ from time $t$ to the final time $T$ by following the schedule $I[t:T]:=(I[t],...,I[T])$.
Similarly, denote the optimal cost-to-go of following the optimal schedule $I^*[t:T]:=(I^*[t],...,I^*[T])$ by $J_t^{I^*[t:T]}(P[t],I[t-1])$, written as $J_t^*(P[t],I[t-1])$ without ambiguity. % due to uniqueness.
Then, one has that
\begin{equation}
J_t^*(P[t],I[t-1]) = \min_{\forall I[t:T]} J_t^{I[t:T]}(P[t],I[t-1])
\end{equation}
for any $P[t]\in\mathcal{P}_t$ and $I[t-1]\in\mathcal{I}_t$.

Considering the cost function \eqref{eq:cost}, the cost-to-go incurred by following schedule $I[t:T]$ is summarized as
\begin{align}
&J_t^{I[t:T]}(P[t],I[t-1])\notag\\
&= \sum_{\tau=t}^T Q(P[\tau],I[\tau]) + \kappa(I[\tau-1],I[\tau])\notag\\
&=Q(P[t],I[t])+\kappa(I[t-1],I[t]) \notag\\
&\qquad\qquad\qquad + J_{t+1}^{I[t+1:T]}(f_{I[t]}(P[t]),I[t])
\label{eq:Jt}
\end{align}
for all $t=1,...,T-1$ and 
\begin{equation}
J_T^{I[T]}(P[T],I[T-1]) = Q(P[T],I[T])+\kappa(I[T-1],I[T]).
\label{eq:JT}
\end{equation}
Define $J_{T+1}^{I[T+1:T]}(f_{I[T]}(P[T]),I[T])=0$, then the above representation \eqref{eq:Jt} holds for all $t=1,...,T-1$ and $t=T$, i.e., \eqref{eq:JT} is a special form of \eqref{eq:Jt}.
Similarly, the optimal cost-to-go from time $t$ can be written as 
\begin{multline}
J_t^*(P[t],I[t-1]) = Q(P[t],I^*[t]) + \kappa(I[t-1],I^*[t]) \\
+ J_{t+1}^*(f_{I^*[t]}(P[t]),I^*[t])
\end{multline}
for all $t=1,...,T$ where $J^*_{T+1}(f_{I^*[T]}(P[T]),I^*[T])=0$.

The Bellman's principle of optimality gives
\begin{multline}
J_t^*(P[t],I[t-1]) = \min_{\forall I[t]} \Big( Q(P[t],I[t]) + \kappa(I[t-1],I[t]) \\
+ J_{t+1}^*(f_{I[t]}(P[t]),I[t])\Big)
\label{eq:bmo}
\end{multline}
for any $P[t]\in\mathcal{P}_t$ and $I[t-1]\in\mathcal{I}_t$.
Then, the optimal schedule $I^*[t]$ is given by
\begin{multline}
I^*_t(P[t],I[t-1]) = \arg\min_{\forall I[t]}\Big(Q(P[t],I[t]) + \kappa(I[t-1],I[t])\\
+ J_{t+1}^*(f_{I[t]}(P[t]),I[t]) \Big).
\label{eq:bma}
\end{multline}
The above \eqref{eq:bma} implies that if one has calculated $J_t^*(P[t],I[t-1])$, then the optimal switching schedule $I^*[t]=I^*_t(P[t],I[t-1])$ can be obtained readily for any $P[t]\in\mathcal{P}_t$, $I[t-1]\in\mathcal{I}_t$, and $t=1,...,T$.

\subsection{Closed-Loop Hierarchical Operation}
Since it is difficult to solve the optimal cost-to-go $J_t^*(P[t],I[t-1])$ in closed form exactly from \eqref{eq:bmo}, approximation approach is employed in the following.
Use a neural network to approximate the optimal cost-to-go as follows
\begin{multline}
\tilde{J}_t^*(P[t],I[t-1]) = (W_{t}^{I_{t-1}})^T\phi(P[t])\approx J_t^*(P[t],I[t-1]),\\
\quad\forall t=1,...,T,\quad\forall I_t\in\mathcal{I}_t, \quad\forall P[t]\in\mathcal{P}_t
\label{eq:apprx}
\end{multline}
where $\phi(\cdot)=[\phi_1(\cdot),...,\phi_M(\cdot)]^T$ with each element $\phi_{i}(\cdot):\mathbb{R}^{N+2}\to\mathbb{R}$ being a continuous basis function such as polynomials.
$(W_{t}^{I_{t-1}})^T\in\mathbb{R}^M$ is the weight matrix, and $M$ is the number of neurons.
With approximated optimal cost-to-go, the optimal schedule $I^*[t]$ can be obtained immediately by
\begin{multline}
\tilde{I}^*_t(P[t],I[t-1]) := \arg\min_{\forall I[t]}\Big(Q(P[t],I[t]) + \kappa(I[t-1],I[t])\\
+ \tilde{J}_{t+1}^*(f_{I[t]}(P[t]),I[t]) \Big).
\end{multline}

Then we have the following CLHO algorithm summarized in Algorithm \ref{alg:CLHO}. 
%In approximation \eqref{eq:apprx}, neural network with continuous active functions has been used to approximate the optimal cost-to-go.
The following Theorem \ref{thm:cont} shows that the optimal cost-to-go is continuous w.r.t. $P[t]$. Hence, the neural network approximation in \eqref{eq:apprx} is possible according to Weierstrass's theorem and \cite{Hornik:1989um}. %, otherwise, the approximation is not guaranteed to be uniform. 
\begin{algorithm}[htbp]
\small
\caption{Closed-Loop Hierarchical Operation (CLHO) Algorithm}
\SetAlgoLined
\SetKwProg{Train}{Training Phase}{}{End}
\SetKwProg{Schedule}{Scheduling Phase}{}{End}

\KwData{Renewable power and load forecasting model \eqref{eq:rpf} \eqref{eq:lf};
	Initial status of $P[0]$ and $I[0]$.}
\KwResult{Optimal schedule $(I^*[t])_{t=1}^T$ of hybrid system \eqref{eq:hybrid} w.r.t. objective \eqref{eq:cost}.}

For $\tau=1,...,T$, if $\tilde{J}_{\tau+1}^*(\cdot,\cdot)$ has been trained, go to \textit{Scheduling Phase}, otherwise go to \textit{Training Phase}.
\Train{}{
Step 1: For $t=T,...,\tau+1$ repeat Step 2 and Step 3\;
Step 2: Select a set of $P[t]\in\mathcal{P}_t$ randomly\;
Step 3: Tune $W_{t}^{I_{t-1}}$ by using the method of least squares such that
	\begin{multline}
	\tilde{J}_t^*(P[t],I[t-1]) := (W_{t}^{I_{t-1}})^T\phi(P[t])\\
	\approx \min_{\forall I[t]}\Big(Q(P[t],I[t]) + \kappa(I[t-1],I[t])\\
	+ \tilde{J}_{t+1}^*(f_{I[t]}(P[t]),I[t]) \Big).\label{eq:train}
	\end{multline}}
\Schedule{}{
Step 4: Input $P[\tau]$ and $I[\tau-1]$ of current time $\tau$\;
Step 5: Calculate the optimal schedule $I^*[t]$ via
\begin{multline}
I^*[\tau] := \tilde{I}^*_t(P[\tau],I[\tau-1]) \\
= \arg\min_{\forall I[\tau]}\Big(Q(P[\tau],I[\tau]) + \kappa(I[\tau-1],I[\tau])\\
+ \tilde{J}_{\tau+1}^*(f_{I[\tau]}(P[\tau]),I[\tau]) \Big).\label{eq:schedule}
\end{multline}}

\label{alg:CLHO}
\end{algorithm}

\begin{thm}
The optimal cost-to-go $J_t^*(P[t],I[t-1])$ of hybrid system \eqref{eq:hybrid} minimizing cost function \eqref{eq:cost} is continuous w.r.t. $P[t]$ for any $I[t-1]\in\mathcal{I}_{t-1}$ and $t=1,...,T$.
\label{thm:cont}
\end{thm}
To prove Theorem \ref{thm:cont}, we need to prove the following lemma first.
\begin{lem}
The dynamic $P[t]=f_I(P[t-1],P^M_{\mDG}[t],P_{\mD}[t])$ of the hybrid system \eqref{eq:hybrid} is continuous w.r.t. $P[t-1]$ for any feasible schedule $I\in\mathcal{I}_t$.
\label{lm:cont}
\end{lem}
\begin{IEEEproof}[Proof of Lemma \ref{lm:cont}]
%\begin{comment}
Note that if we can obtain the solution of optimization problem \eqref{eq:lop} in explicit form for any $I$, then the continuity of dynamic $f_I(\cdot,P^M_{\mDG}[t],P_{\mD}[t]): \mathbb{R}^{N+2}\mapsto\mathbb{R}^{N+2}$ can be easily obtained by analyzing the explicit solution.
However, it is difficult to obtain such an explicit solution.
Thus, we turn to the conditions on which the solution (not necessarily in explicit form) can be guaranteed.

%Since it is difficult to obtain the solution of \eqref{eq:lop} in explicit form, we turn to the conditions on which the solution can be guaranteed.
%\end{comment}
For any feasible $I\in\mathcal{I}_t$, function $Q(P[t],I)$ is convex w.r.t. $P[t]$ since $C_{\mF,n}(P_n[t])$, $E_n(P_n[t])$, $C_{\mDG}(P_{\mDG}[t])$, and $C_{\mDR}(P_{\mDR}[t])$ are convex for all $n$.
The power balance constraint \eqref{eq:pbc} is linear, and all the other constraints \eqref{eq:sr1}--\eqref{eq:drc} are affine.
Thus, the optimization problem \eqref{eq:lop} has strong duality.
In addition, the function $Q(P[t],I)$ is differentiable w.r.t. $P[t]$ for any feasible $I\in\mathcal{I}_t$.
The convex optimization theory \cite{Boyd:2014uz} gives that a solution of \eqref{eq:lop} is optimal iff the solution satisfies the Karush-Kuhn-Tucker (KKT) conditions.
%\cite{Boyd:2014uz}的page 244和\cite{Bertsekas:2018wk}的page 5

Write the optimization problem \eqref{eq:lop} in general form as
\begin{align}
\min\ & Q(P[t],I[t])\\
\st\ & h(P[t])=0,\ g(P[t])\leq0,
\end{align}
where $h(P[t]):=\sum_{n\in\mathcal{C}} P_n[t] + P_{\mDG}[t] + P_{\mDR}[t] - P_{\mD}[t]$, $g(P[t])$ is the stack vector of functions $\sum_{n\in\mathcal{C}} P_n^m + P_{\mDG}[t] + P_{\mDR}[t] - P_{\mD}[t]+ R^m[t]$, $P_{\mD}[t]+ R^M[t]-\sum_{n\in\mathcal{C}} P_n^M - P_{\mDG}[t] - P_{\mDR}[t]$, $P^m_n - P_n[t]$, $P_n[t] - P^M_n$, $-P_{\mDG}[t]$, $-P_{\mDR}[t]$, $P_{\mDG}[t]-P^M_{\mDG}[t]$, $P_{\mDR}[t]-P^M_{\mDR}[t]$, and $(1-\eta_{\max})P_{\mDG}[t]-\eta_{\max} \sum_{n\in\mathcal{C}[t]}P_n[t]$ for $\forall n\in \mathcal{C}[t]:=\{n|I_n[t]=1\text{ or }P_n[t]>0\}$, $-R^{\mD}_n-P_n[t]+P_n[t-1]$ and $P_n[t]-P_n[t-1]-R^{\mU}_n$ for $\forall n\in \mathcal{C}[t]\cap\mathcal{C}[t-1]$.
The KKT conditions are as follows
\begin{multline}
\nabla Q(P[t],I[t]) + \langle \lambda,\nabla h(P[t]) \rangle + \langle \mu, \nabla g(P[t]) \rangle =0,\\
h(P[t])=0,\ \mu\geq 0,\ g(P[t])\leq 0,\ \langle \mu,g(P[t]) \rangle =0.
\label{eq:kkt}
\end{multline}
Define $(P_t^*,\lambda^*,\mu^*)$ the solution of \eqref{eq:kkt} with $P_t^*:=P^*[t]$, which is the optimal primal-dual solution of problem \eqref{eq:lop} as argued previously.
Denote the solution of \eqref{eq:kkt} as $(\bar{P}_t,\bar{\lambda},\bar{\mu})$ when $g(P[t])$ is perturbed by some $\mathbf{c}$ such that $-R^{\mD}_n-P_n[t]+P_n[t-1]\leq c$ and $P_n[t]-P_n[t-1]-R^{\mU}_n\leq -c$ for $\forall n\in \mathcal{C}[t]\cap\mathcal{C}[t-1]$ while other inequalities remains unchanged.
That is to say, $(\bar{P}_t,\bar{\lambda},\bar{\mu})$ is the solution of perturbed KKT conditions
\begin{multline}
\nabla Q(P[t],I[t]) + \langle \lambda,\nabla h(P[t]) \rangle + \langle \mu, \nabla g(P[t]) \rangle =0,\\
h(P[t])=0,\ \mu\geq 0,\ g(P[t])\leq \mathbf{c},\ \langle \mu,g(P[t]-\mathbf{c}) \rangle =0.
\label{eq:pkkt}
\end{multline}
The upper Lipschitz stability \cite{Izmailov:2012di} gives that there is a neighborhood of $(P^*_t,\lambda^*,\mu^*)$ such that for any $\mathbf{c}$ close enough to $\mathbf{0}$, it holds
\begin{equation}
\|\bar{P}_t - P^*_t\|+\text{dist}((\lambda^*,\mu^*),(\bar{\lambda},\bar{\mu})) \leq \|\mathbf{c}\|.
\end{equation}
Therefore, for any $\epsilon>0$, there exists $0<\delta<\epsilon$, if $\|P_n[t-1]-(P_n[t-1]-c)\|=c<\delta$, $\forall n\in \mathcal{C}[t]\cap\mathcal{C}[t-1]$, then
$\|\bar{P}_t-P^*_t\|\leq\|\mathbf{c}\|<\delta<\epsilon$. Thus, one obtains the continuity of $f_I(P[t-1],P^M_{\mDG}[t],P_{\mD}[t])$ w.r.t. $P[t-1]$.
\end{IEEEproof}
Now we are ready to prove the Theorem \ref{thm:cont}.
\begin{IEEEproof}[Proof of Theorem \ref{thm:cont}]
Firstly, one has the continuity of $J^*_T(P[T],I[T-1])$ w.r.t. $P[T]$ since $Q(P[T],I^*[T])$ is continuous and $\kappa(I[T-1],I^*[T])$ is constant w.r.t. $P[T]$.
Secondly, it can be shown by mathematical induction that $J^*_t(P[t],I[t-1])$ is continuous w.r.t. $P[t]$ for any $I[t-1]\in\mathcal{I}_{t-1}$ and $t=1,...,T-1$.

Assume that $J^*_{t+1}(P[t+1],I[t])$ is continuous w.r.t. $P[t+1]$.
According to Bellman's principle of optimality, one has that
\begin{multline}
J_t^*(P[t],I[t-1]) = \min_{\forall I[t]} \Big( Q(P[t],I[t]) + \kappa(I[t-1],I[t]) \\
+ J_{t+1}^*(f_{I[t]}(P[t]),I[t])\Big),
\end{multline}
which is continuous since $F(P[t],I[t-1],I[t]):=Q(P[t],I[t]) + \kappa(I[t-1],I[t]) + J_{t+1}^*(f_{I[t]}(P[t]),I[t])$ is continuous for any $I[t-1]\in\mathcal{I}_{t-1}$ and $I[t]\in\mathcal{I}_t$ where $|\mathcal{I}_t|$ is finite.
The reason is as follows.
Assume that $|\mathcal{I}_t|=2$, then it can be verified that $\min\{F_1,F_2\}=\frac{1}{2}(F_1+F_2-|F_1-F_2|)$ is continuous where $F(P[t],I[t-1],I[t])$ for $I[t]\in\mathcal{I}_t$ are denoted by $F_1$ and $F_2$ respectively.
By mathematical induction, one has that $\min_{i=1,...,|\mathcal{I}_t|} \{F_i\}$ is also continuous given the continuity of $F_i$ with $i$ in finite set $|\mathcal{I}_t|$.
Therefore, $J_t^*(P[t],I[t-1])$ is continuous w.r.t. $P[t]$ for any $I[t-1]\in\mathcal{I}_{t-1}$ given the continuity of $J_{t+1}^*(P[t+1],I[t])$.
This completes the proof.
\end{IEEEproof}

The above Theorem \ref{thm:cont} shows that the optimal cost-to-go $J_t^*(P[t],I[t-1])$ is continuous w.r.t. $P[t]$.
Hence, $\tilde{J}_t^*(P[t],I[t-1])$ in Algorithm \ref{alg:CLHO} uniformly approximates $J_t^*(P[t],I[t-1])$ by using training methods such as least squares according to Weierstrass's theorem and \cite{Hornik:1989um}.
That is to say, the CLHO algorithm provides a means of computing the optimal cost-to-go function $J_t^*(P[t],I[t-1])$. % which can not be solved in closed form exactly from \eqref{eq:bmo}.
Note that the neural networks trained in Algorithm \ref{alg:CLHO} produce the representation of \textit{function} $\tilde{J}_t^*(P[t],I[t-1])$ (rather than \textit{numbers}). %, the optimal schedule is calculated by \eqref{eq:schedule} for any $P[t]$ and $I[t-1]$, $t=1,...,T$.
The representation of function $\tilde{J}_t^*(P[t],I[t-1])$ implies that for any (disturbed) $P[t]$ and $I[t-1]$, the schedule produced by \eqref{eq:schedule} without retraining is still optimal for the following $[t,T]$.
Therefore, we say Algorithm \ref{alg:CLHO} produces optimal solution in closed-loop style.

\section{Simulation Cases\label{sec:Sim}}

In this section, simulation cases are presented to study the behavior and effectiveness of the proposed CLHO algorithm.
\subsection{Example 1}
\begin{figure}[htbp]
\centering
\includegraphics[width=88mm,keepaspectratio]{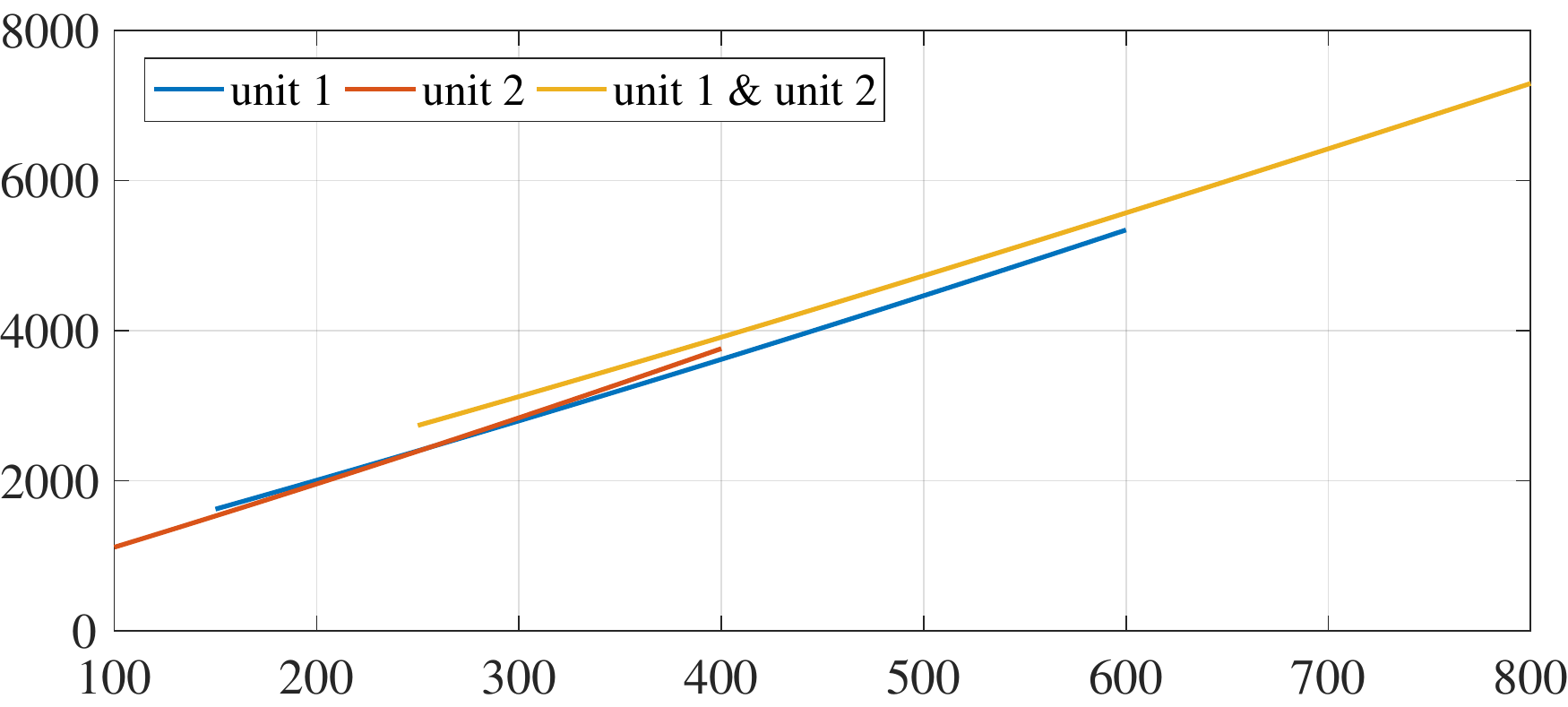}%88mm
\caption{Generation cost for serving demand 100 -- 800 MW.}
\label{plt:cost}
\end{figure}

In this simulation example, the closed-loop behavior of the CLHO algorithm is studied.
A relatively simple UCD problem is selected in this example where the optimal solution can be easily spotted or found with method of exhaustion for ease of studying.
Specifically, this example contains only two thermal units while the ramp limit and spinning reserves are not considered.
The parameters of the thermal units are shown in Table \ref{tab:P}.
The demand is shown in table \ref{tab:demand}.
The generation costs of the thermal units for serving the demand varying from 100 MW to 800 MW are plotted in Fig. \ref{plt:cost}.
 
\begin{table}[htbp]
\centering\small
	\caption{Parameters of thermal units for example 1}
\begin{tabular}[center]{c c c c c c p{.38cm} p{.38cm} p{.38cm}}
	\toprule[1pt]\midrule[0.3pt]
	Unit & $a_n$ & $b_n$ & $c_n$ & $P_n^m$ & $P_n^M$ & $C_{\mb,n}$ & $C_{\mf,n}$ & $C_{D,n}$\\ \midrule
	
	1 & 0.00142 & 7.20 & 510 & 150 & 600 & \ 0 & \ 0 & \ 0\\
	2 & 0.00194 & 7.85 & 310 & 100 & 400 & \ 0 & \ 0 & \ 0\\
	
	\midrule[0.3pt]\bottomrule[1pt]
\end{tabular}
\label{tab:P}
\end{table}

\begin{table}[htbp]
\centering\small
	\caption{Forecasted power demand for example 1}
\begin{tabular}[center]{c c c c c c c c}
	\toprule[1pt]\midrule[0.3pt]
	Time (h) & 0 & 1 & 2 & 3 & 4 & 5 & 6\\ \midrule
	
	Demand (MW) & 200 & 200 & 350 & 350 & 700 & 700 & 700\\
	
	\midrule[0.3pt]\bottomrule[1pt]
\end{tabular}
\label{tab:demand}
\end{table}

\begin{figure}[htbp]
\centering
\subfloat[][Only unit 2 is committed]{\includegraphics[width=88mm,keepaspectratio]{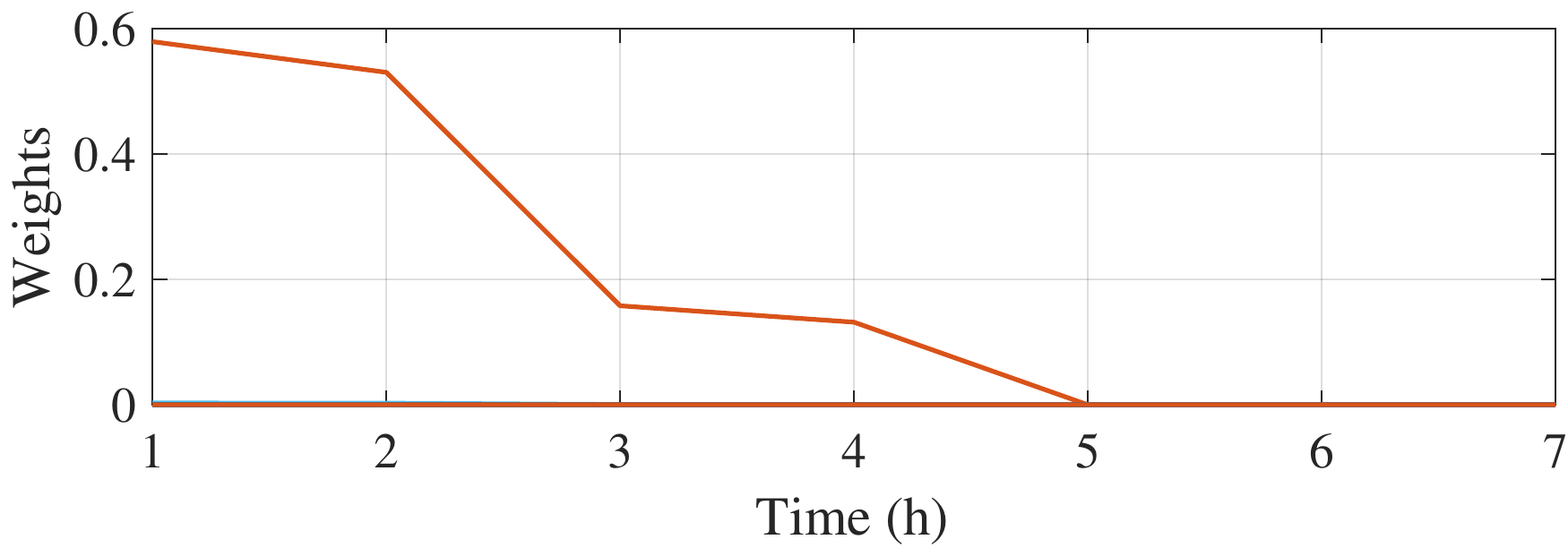}}\\%keepaspectratio
\subfloat[][Only unit 1 is committed]{\includegraphics[width=88mm,keepaspectratio]{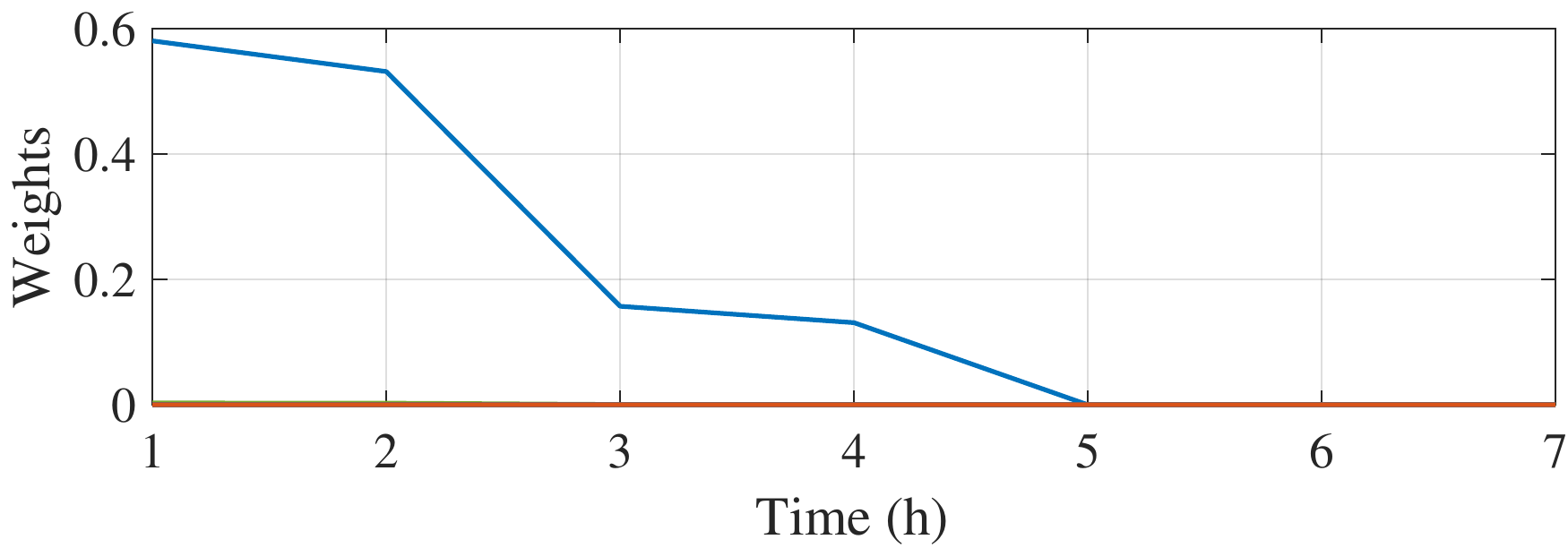}}\\
\subfloat[][Both units are committed]{\includegraphics[width=88mm,keepaspectratio]{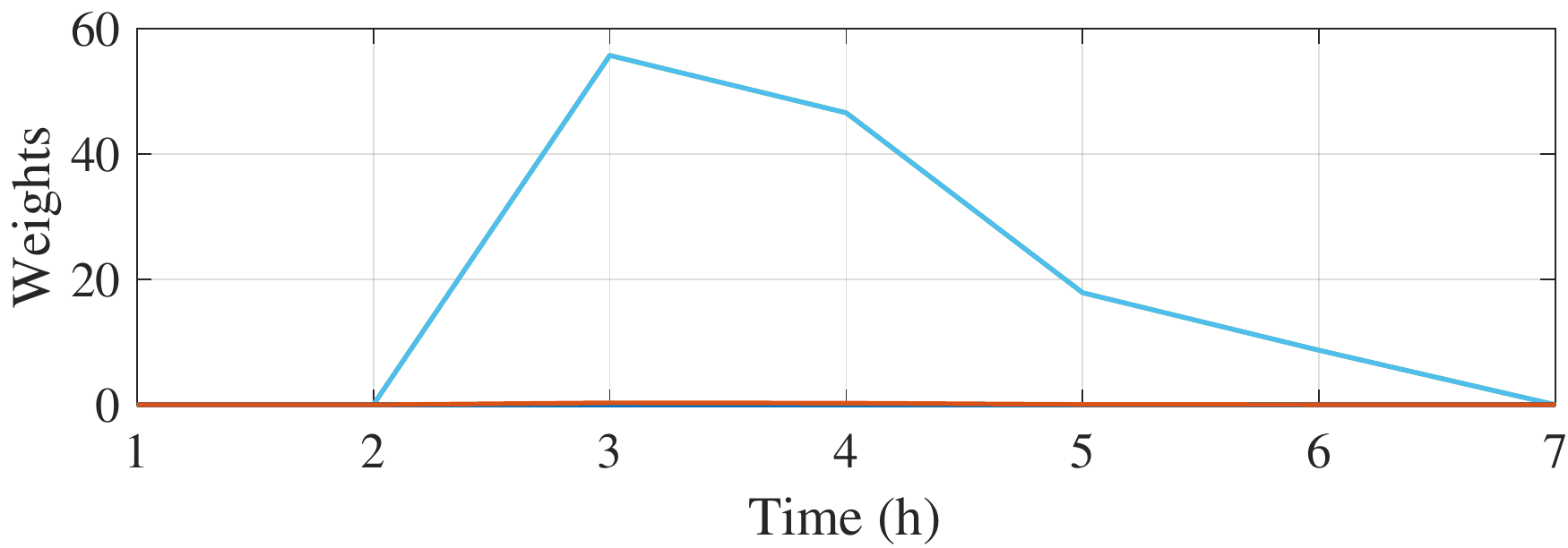}}
\caption{Weight histories of the trained neural networks.\label{plt:W}}
\end{figure}

\textit{Training phase:}
The basis functions are selected as $\phi(P[t])=[P_1^2[t], P_2^2[t], P_1[t], P_2[t], 1]^T$.
At each time step $t$ and schedule $I_{t-1}$, one hundred possible power output combinations are sampled randomly.
The weights $W_t^{I_{t-1}}$ of the neural networks are decided by using least squares.
The resulting weights of the neural networks are shown in Fig. \ref{plt:W}.
%CVX package is employed to solve the lower level optimization problem \eqref{eq:lop}.
The training phase takes 42 s in MATLAB 2016a on a PC with Intel Core i3-4130 CPU and 8 GB of RAM.

\textit{Scheduling phase:}
Once trained, the neural networks are ready to be used to schedule the generation units. In the following, four simulation cases are presented to study the closed-loop behavior of the CLHO algorithm.

\subsubsection{Case 1}
In the first case, the initial conditions are selected as $P[0]=[0, 200]^T$ MW. The scheduling results are shown in Fig. \ref{plt:init1}.
Since the start-up and shutdown costs are zero as shown in Table \ref{tab:P}, one has that the optimal scheduling at each time step is the optimal power output combinations for that time step.
According to Fig. \ref{plt:cost}, the optimal schedule for demand of 200 MW, 350 MW, and 700 MW are $I=[0, 1]^T$, $I=[1, 0]^T$, and $I=[1, 1]^T$ respectively.
It can also be verified that the corresponding optimal generation outputs are $P=[0, 200]^T$ MW, $P=[350, 0]^T$ MW, and $[500.9, 199.1]^T$ MW.
Therefore, the scheduling results as shown in Fig. \ref{plt:init1} is optimal with total generation cost \$\,23168.3.

\subsubsection{Case 2}
In the second case, the initial conditions remain the same as $P[0]=[0, 200]^T$ MW. However, the power outputs at time $t=2$ are $P[2]=[200, 150]^T$ MW due to external disturbance. No re-train is needed in this case. The resulting history of generation outputs are shown in Fig. \ref{plt:init1dist}.
The dashed lines represent the scheduling results without disturbance while the real lines show the actual power outputs under disturbance.
The results show that though the power outputs are $P[2]=[200, 150]^T$ MW at time $t=2$ due to disturbance, the following power outputs are still optimal comparing with Fig. \ref{plt:init1}.
Thus, the scheduling phase of the CLHO algorithm is closed-loop to external disturbances.

\subsubsection{Case 3}
\begin{figure}[htbp]
\centering
%\minipage{0.30\textwidth}
\includegraphics[width=88mm,keepaspectratio]{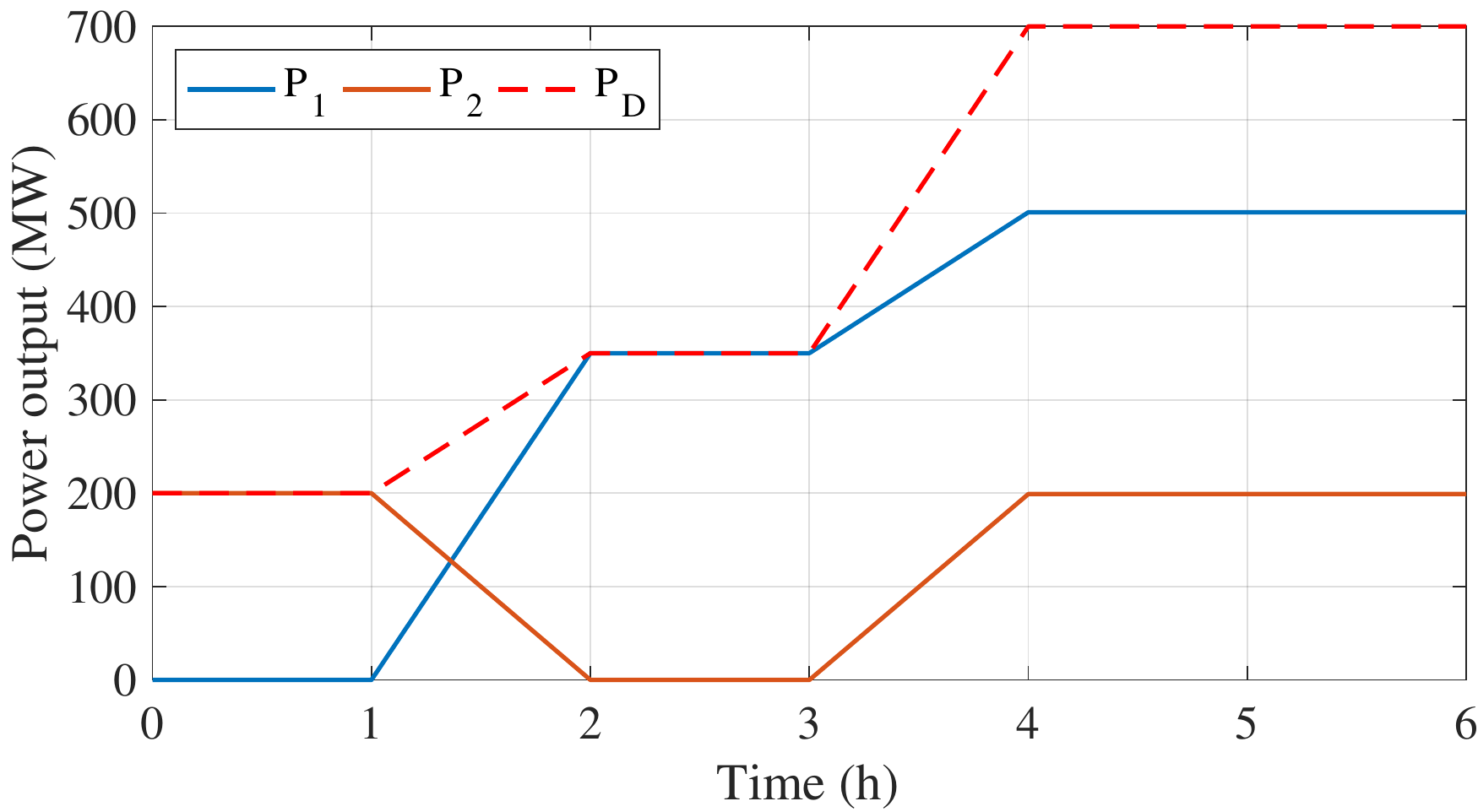}%88mm
\caption{Scheduling results with initial $P[0]=[0, 200]^T$ MW.}
\label{plt:init1}
%\endminipage\hspace{5mm}
%\minipage{0.30\textwidth}
\centering
\includegraphics[width=88mm,keepaspectratio]{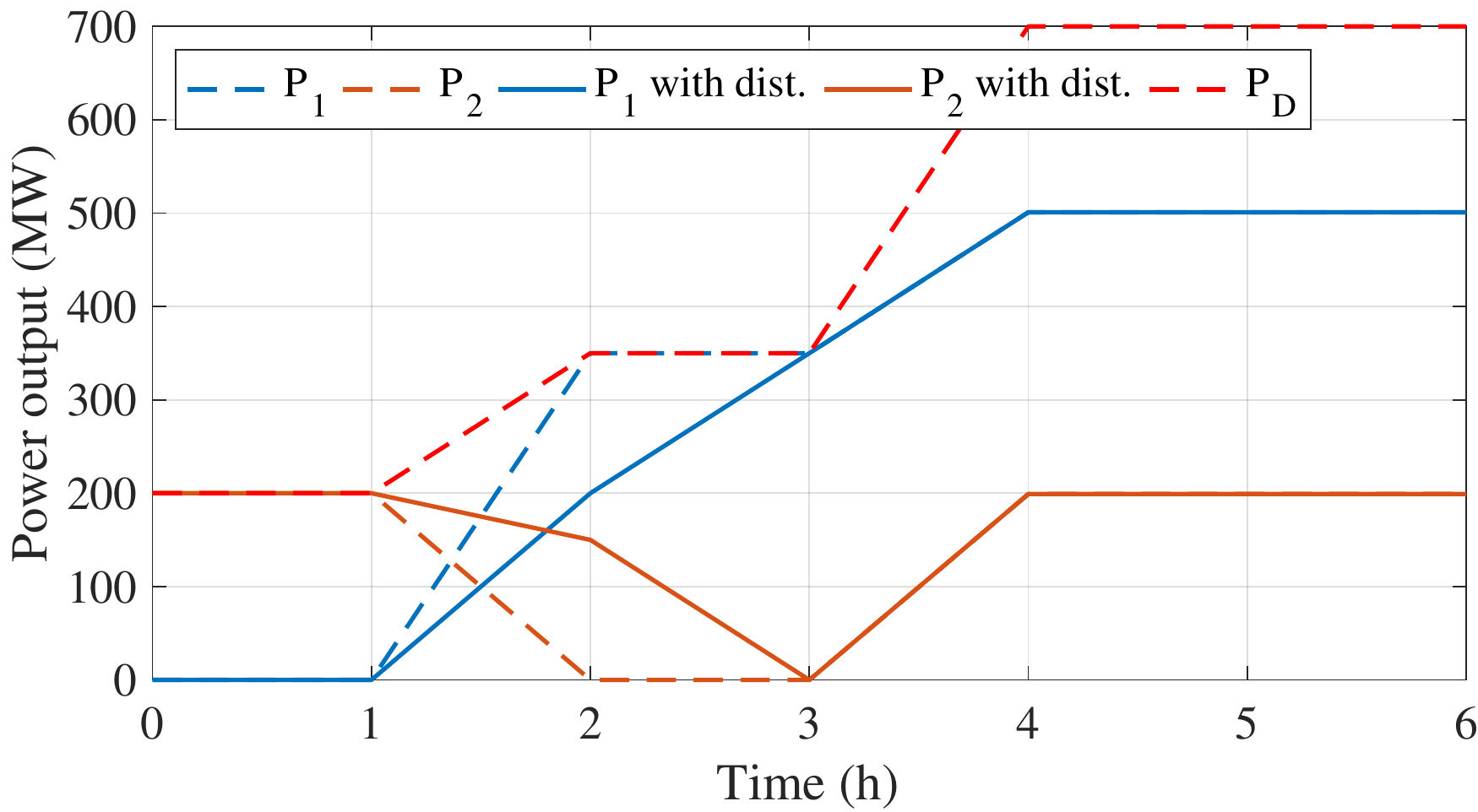}%88mm
\caption{Scheduling results under disturbance at time $t=2$.}
\label{plt:init1dist}
%\endminipage\hspace{5mm}
%\minipage{0.30\textwidth}
\includegraphics[width=88mm,keepaspectratio]{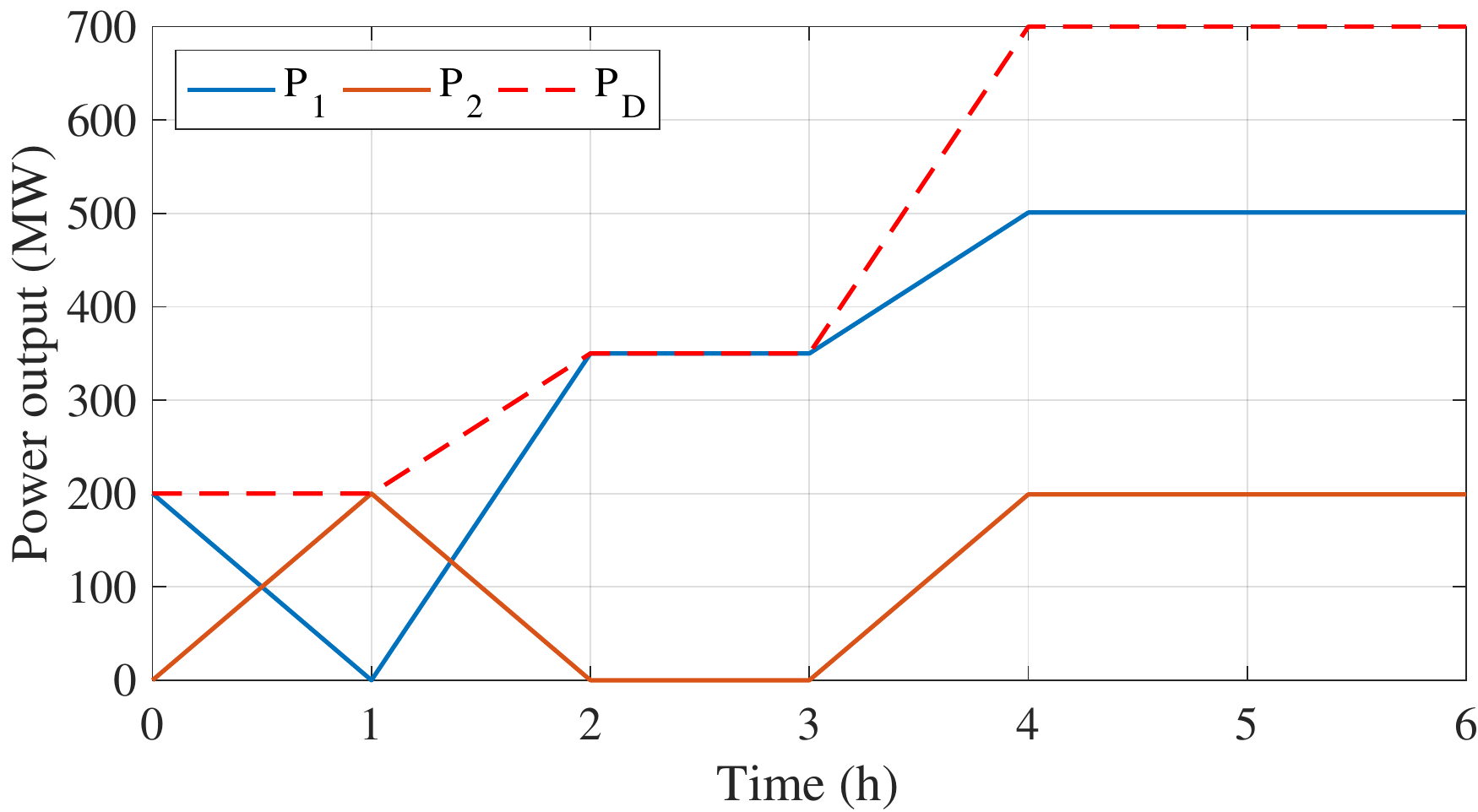}%88mm
\caption{Scheduling results with initial $P[0]=[200, 0]^T$ MW.}
\label{plt:init2}
%\endminipage
\end{figure}

In the third case, the initial conditions are changed to $P[0]=[200, 0]^T$ MW.
No re-train is needed in this case. The scheduling results are shown in Fig. \ref{plt:init2}.
Since there is no start-up and shutdown costs, the generation outputs are scheduled to the optimal one at first time step $t=1$. In what follows, the scheduling results remain optimal comparing with Fig. \ref{plt:init1}.
Thus, the scheduling phase of the CLHO algorithm is also closed-loop to different initial conditions.
%Note that the UCD is usually conducted under the assumption that the power demand and distributed energy sources are accurately forecasted.
%In prevalent real-time operation, the system operator also rely on ancillary services market and real-time market to cope with the forecasting errors.

\subsubsection{Case 4}
\begin{figure}[htbp]
\centering
\includegraphics[width=88mm,keepaspectratio]{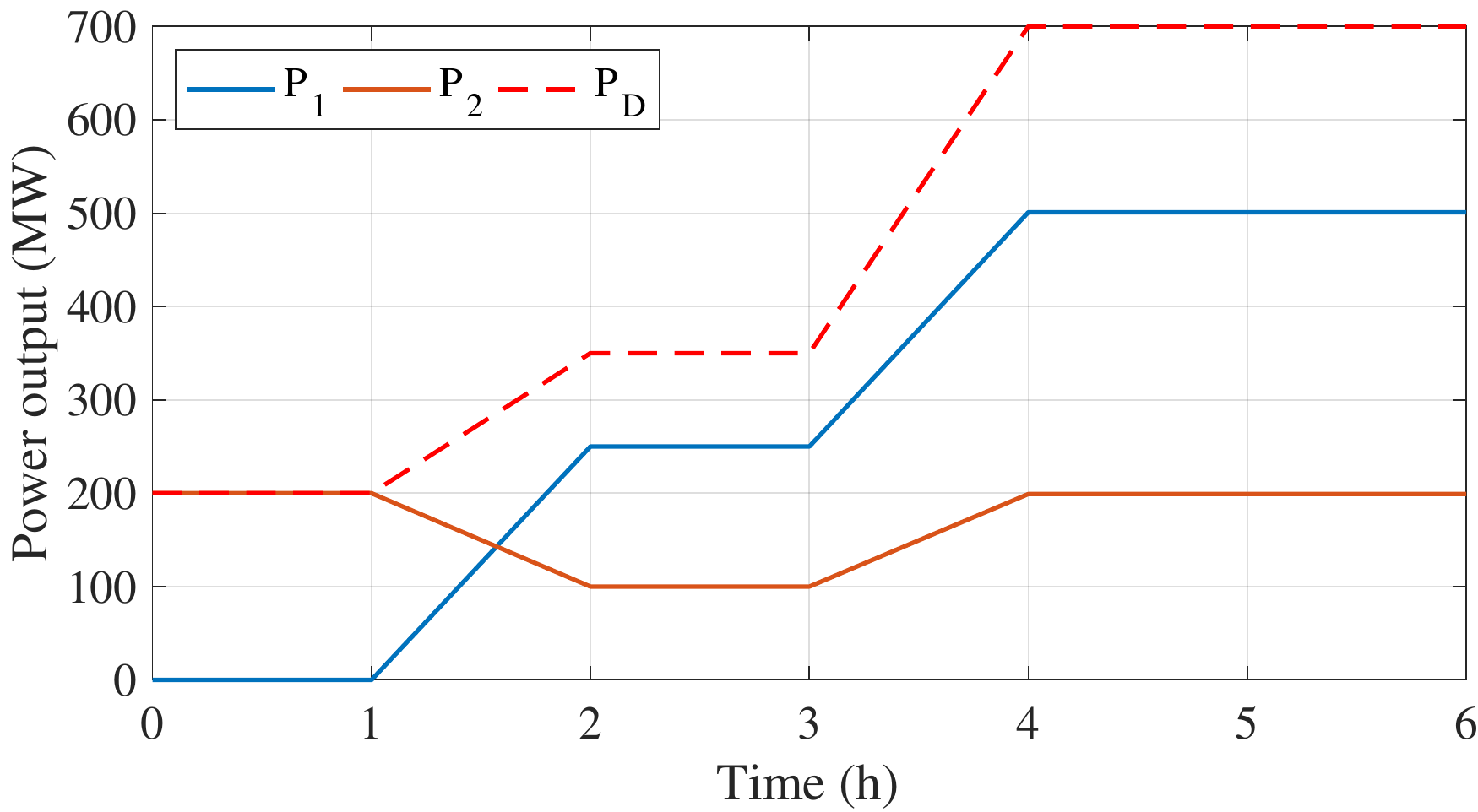}%88mm
\caption{Scheduling results considering start-up/shutdown costs.}
\label{plt:init1swcost}
\end{figure}

In the fourth case, the start-up costs are set as $C_{\mU,1}(\tau_1)=300 \cdot \tau_1$ and $C_{\mU,2}(\tau_2)=200 \cdot \tau_2$.
The shutdown costs are set as $C_{\mD,1}=600$ and $C_{\mD,2}=400$.
Due to the changes of switching cost in \eqref{eq:train} incurred by the start-up and shutdown costs, the neural networks are re-trained in this case.
After training the neural networks, the scheduling phase is carried out.
The scheduling results are shown in Fig. \ref{plt:init1swcost}.
Comparing the results with that shown in Fig. \ref{plt:init1}, one has that the unit 2 is not shutdown at $t=2$ and start-up again at $t=3$ due to the shutdown and start-up costs.
The total cost from time $t=0$ to time $t=6$ incurred by this scheduling results is \$\,24386.7, which is less than the total cost \$\,24568.3 incurred by the scheduling results in Fig. \ref{plt:init1} including non-zero start-up and shutdown costs.
In fact, the scheduling results given in Fig. \ref{plt:init1swcost} is optimal, which can be verified in the following with method of exhaustion.
All the possible commitment and the corresponding minimal total generation costs for both the first case and this case are shown in Table \ref{tab:commitment}.
In the table, number 1, 2, and 3 represent the commitment $[0,1]^T$, $[1,0]^T$, and $[1,1]^T$ respectively.
According to the results shown in Table \ref{tab:commitment}, the optimal commitment for $t=1$ to $t=6$ is $133333$ and the minimal total generation cost is \$ 24386.7, which are the same as the scheduling results shown in Fig. \ref{plt:init1swcost} and the corresponding total cost.

Similarly, with the help of Fig. \ref{plt:cost}, it can be easily verified that all the scheduling results for the previous three simulation cases are also optimal under the corresponding set up.
The neural networks are not retrained in the first three cases, revealing an important feature that the proposed algorithm produces scheduling results in closed-loop style.
Since the optimal schedule can be readily produced without retraining, the CLHO algorithm may contribute to the real-time operation of microgrids.
As comparison, most existing methods such as \cite{UCviaPL,UCviaLR,UCviaPSO} require reoptimization for different initial conditions and external disturbances.

Simulation comparison with traditional methods is not presented due to the following fundamental and clear distinctions. 
Firstly, the CLHO algorithm can be used differently from traditional methods.
Namely, the time-consuming training phase can be conducted day-ahead while in the scheduling phase, no retraining is needed for different initial conditions and disturbances as shown above.
Secondly, the CLHO algorithm finds the optimal solution as shown by the theoretical analysis and simulations, while most traditional methods usually produce suboptimal solutions.

\begin{table}[htbp]
\centering\footnotesize
	\caption{Minimal generation costs for different schedule}
\begin{tabular}[center]{p{1cm} p{.75cm} p{.75cm} p{.75cm} p{.75cm} p{.75cm} p{.75cm}}
	\toprule[1pt]\midrule[0.3pt]
	\scriptsize{$I[1:6]$} & 111333 & 112333 & 113333 & 121333 & 122333 & 123333\\
	\scriptsize{Case 1} & 23350.7 & 23259.5 & 23568.7 & 23259.5 & \underline{23168.3} & 23477.5\\
	\scriptsize{Case 4} & 24550.7 & 24759.5 & 24468.7 & 25359.5 & 24568.3 & 24677.5\\\midrule[0.3pt]
	
	\scriptsize{$I[1:6]$} & 131333 & 132333 & 133333 & 211333 & 212333 & 213333 \\
	\scriptsize{Case 1} & 23568.7 & 23477.5 & 23786.7 & 23399.9 & 23308.7 & 23617.9\\
	\scriptsize{Case 4} & 25068.7 & 24677.5 & \underline{24386.7} & 25499.9 & 25708.7 & 25417.9\\\midrule[0.3pt]
	
	 \scriptsize{$I[1:6]$} & 221333 & 222333 & 223333 & 231333 & 232333 & 233333\\
	 \scriptsize{Case 1} & 23308.7 & 23217.5 & 23526.7 & 23617.9 & 23526.7 & 23835.9\\
	 \scriptsize{Case 4} &25308.7 & 24517.5 & 24626.7 & 25417.9 & 25026.7 & 24735.9\\
 
	\midrule[0.3pt]\bottomrule[1pt]
\end{tabular}
\label{tab:commitment}
\end{table}

\subsection{Example 2}
In this simulation example, the IEEE 57-bus system is employed to study the effectiveness of the CLHO algorithm.
The forecasted power demand and maximum power generation of DGs over 24 hours is shown in Table \ref{tab:demand24}.
Table \ref{tab:fuel24} and \ref{tab:emission24} show the fuel cost coefficients, generation capacities, start-up and shutdown cost, and emission coefficients of the thermal units.
The parameters of the aggregator managing the DGs and virtual unit of DR are also shown in Table \ref{tab:fuel24}.
The penetration rate $\eta_t$ of the DGs is bounded by $\eta_{\max}=0.05$.
Maximum demand curtailment of the DR is set as 40MW for peak hours ($t=11$\,h and $t=12$\,h) and 10MW for non-peak hours.
To cope with the renewable energy and demand forecast errors, spinning reserves are also considered.
Specifically, the lower and upper spinning reserves are set as 5\% of the forecasted demand $P_{\mD}[t]$ for each hour, namely, $R^{m}[t]=R^{M}[t]=0.05\cdot P_{\mD}[t]$. 
At the first hour $t=1$, the power outputs of unit 1 and unit 2 are set as 500 MW and 200 MW respectively, while assuming zero power outputs for other units and the demand curtailment.

\subsubsection{Case 1}
In the first case, the quotas $Q_n$ is set as 0 ton for each thermal unit, thus the operators have to buy all the emission allowance from the market at a price assumed as $p_e=1$ \$/ton.
The scheduling results are shown in Fig. \ref{plt:init24Q0P1}.
The total emission is 31086 ton.

\subsubsection{Case 2}
In the second case, the quotas $Q_n$ is still set as 0 ton for each thermal unit while the price of emission quotas %in the CET market 
is increased to $p_e=10$ \$/ton.
The scheduling results are shown in Fig. \ref{plt:init24Q0P10}.
Numerical simulation shows that the total emission is 28775 ton, which is less than 31086 ton in the first case.

Since the higher price to purchase emission quota results in a higher emission cost, the thermal units are less preferable than in the first case.
In consequence, the power output of DGs and curtailment of power demand are larger than in the first case, which reduces the total emission.

\subsubsection{Case 3}
\begin{figure}
%\minipage{0.30\textwidth}
\includegraphics[width=88mm,keepaspectratio]{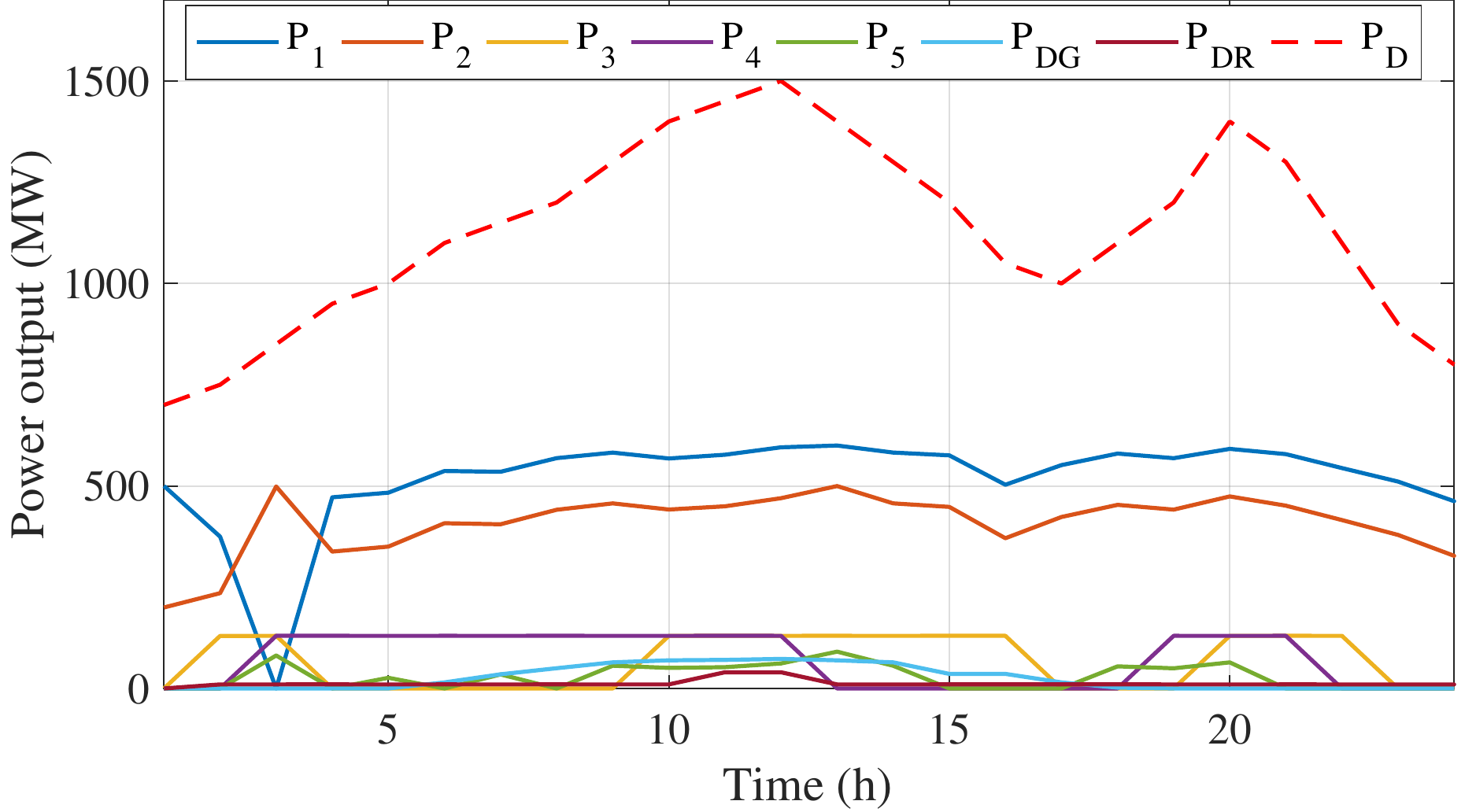}%width=88mm,keepaspectratio
\caption{Scheduling results for case 1 of example 2.}
\label{plt:init24Q0P1}
%\endminipage\hspace{5mm}
%\minipage{0.30\textwidth}
\includegraphics[width=88mm,keepaspectratio]{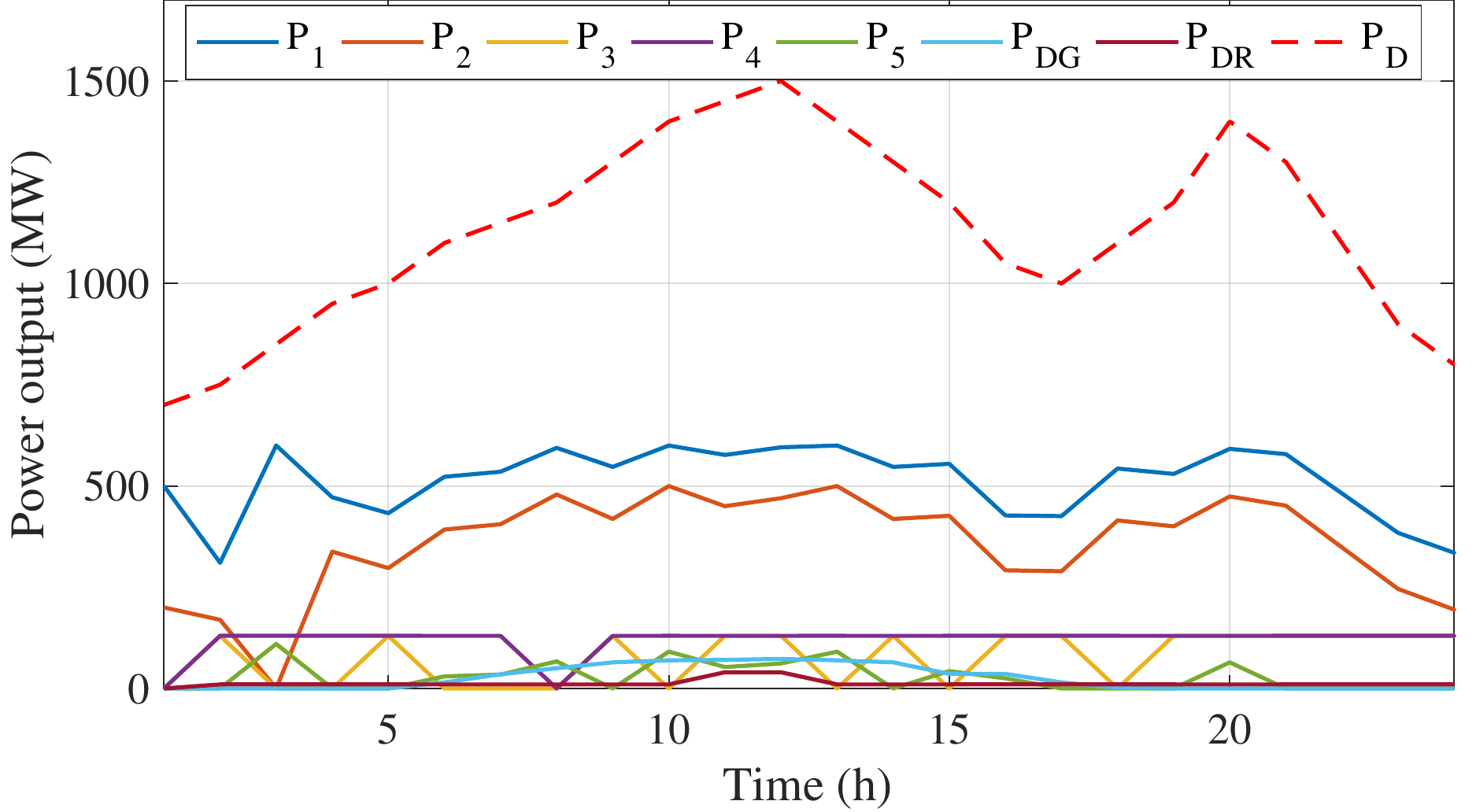}%88mm
\caption{Scheduling results for case 2 of example 2.}
\label{plt:init24Q0P10}
%\endminipage\hspace{5mm}
%\minipage{0.30\textwidth}
\includegraphics[width=88mm,keepaspectratio]{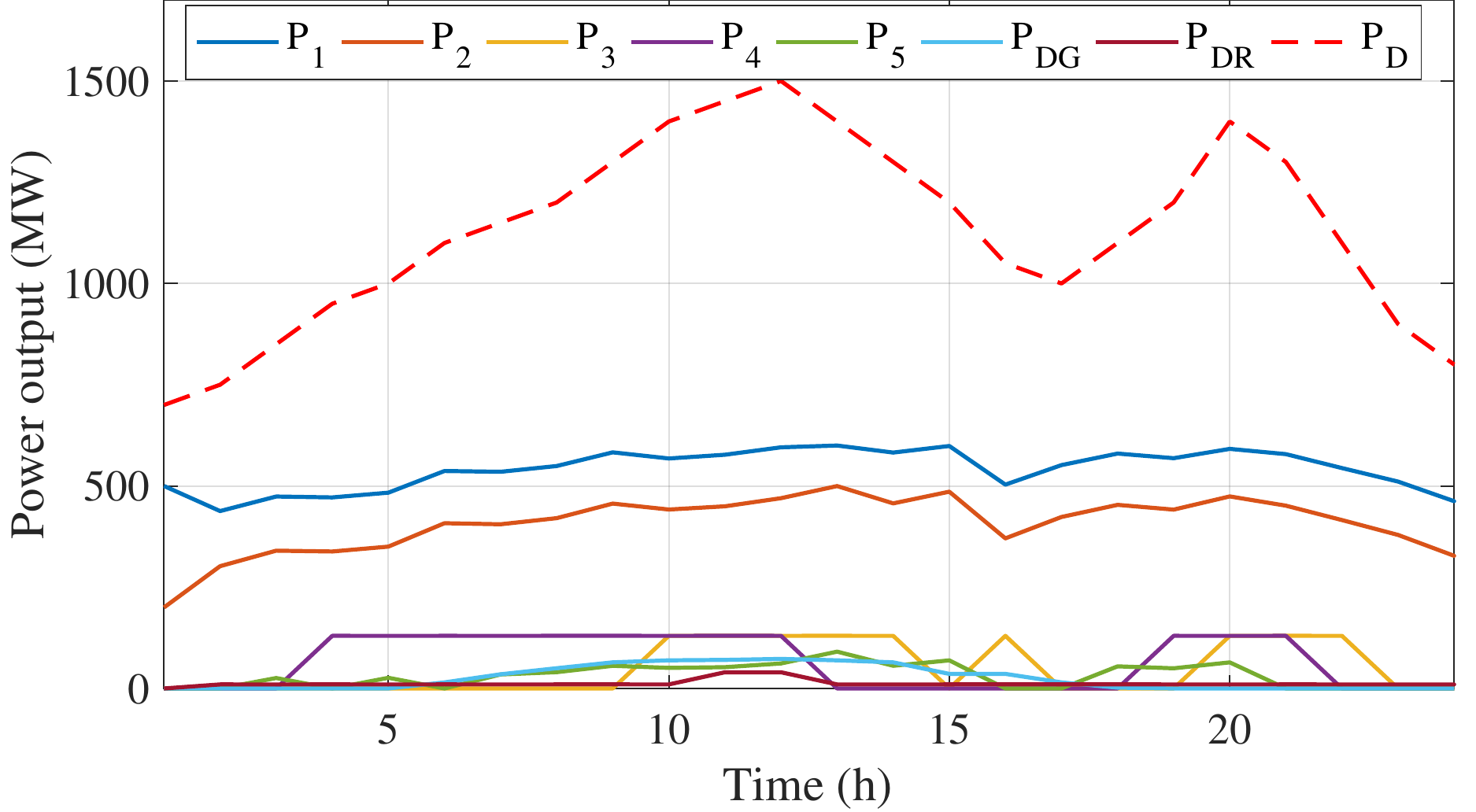}%88mm
\caption{Scheduling results for case 3 of example 2.}
\label{plt:init24Q9P1}
%\endminipage
\end{figure}

In the third case, the quotas $Q_n$ is set as 90\% of the purchased quotas in the first case.
%That is to say, if one use the scheduling results given in Fig. \ref{plt:init24Q0P1}, only 10\% of the allowance is needed to buy from the market compared with 100\% in the first case.
The market price is set as $p_e=1$ \$/ton.
Fig. \ref{plt:init24Q9P1} shows the corresponding scheduling results.
Simulation shows that the total emission is 31404 ton, which is more than the emission in the first case.

Note that the quotas in the first case are less than the quotas here.
The less quotas result in higher emission costs due to the need of buying extra emission allowance.
To reduce the total generation cost, the thermal units are less used in the first case which in turn causes a reduction on the carbon emissions.

In summary, the simulation results show that both the low emission quotas and the high emission trading price help reduce the total amount of carbon emissions.

\begin{table}[htbp]
\centering\footnotesize
	\caption{Forecasted demand and maximum DGs}
\begin{tabular}[center]{c p{.38cm} p{.38cm} p{.38cm} p{.38cm} p{.38cm} p{.38cm} p{.38cm} p{.38cm}}
	\toprule[1pt]\midrule[0.3pt]
	Time (h)& 1 & 2 & 3 & 4 & 5 & 6 & 7 & 8\\
	Demand (MW) & 700 & 750 & 850 & 950 & 1000 & 1100 & 1150 & 1200\\
	Max. DGs (MW) & 0 & 0 & 0 & 0 & 0 & 15 & 35 & 50\\\midrule
	
	Time (h) & 9 & 10 & 11 & 12 & 13 & 14 & 15 & 16\\
	Demand (MW) & 1300 & 1400 & 1450 & 1500 & 1400 & 1300 & 1200 & 1050\\
	Max. DGs (MW) & 72 & 84 & 80 & 88 & 85 & 76 & 36 & 36\\\midrule
	
	Time (h) & 17 & 18 & 19 & 20 & 21 & 22 & 23 & 24\\
	Demand (MW) & 1000 & 1100 & 1200 & 1400 & 1300 & 1100 & 900 & 800\\
	Max. DGs (MW) & 15 & 2 & 0 & 0 & 0 & 0 & 0 & 0\\
	\midrule[0.3pt]\bottomrule[1pt]
\end{tabular}
\label{tab:demand24}
\end{table}

\begin{table}[htbp]
\centering\footnotesize
	\caption{Parameters of (virtual) units for example 2}
\begin{tabular}[center]{c c c c c c p{.38cm} p{.38cm} p{.38cm} c c c}
	\toprule[1pt]\midrule[0.3pt]
	Unit & $a_{(.)}$ & $b_{(.)}$ & $c_{(.)}$ & $P_{(.)}^m$ & $P_{(.)}^M$ & $C_{\mb,n}$ & $C_{\mf,n}$ & $C_{D,n}$\\ \midrule
	1 & 0.00048 & 16.19 & 990 & 150 & 600 & 550 & 100 & 90\\
	2 & 0.00031 & 17.26 & 970 & 100 & 500 & 570 & 110 & 100\\
	3 & 0.00200 & 16.60 & 700 & 20  & 130 & 110 & 85  & 75\\
	4 & 0.00211 & 16.50 & 680 & 20  & 130 & 120 & 95  & 85\\
	5 & 0.00398 & 19.70 & 450 & 25  & 162 & 150 & 100 & 90\\
	DG& 0.01000 & 2.60 & 10 & 0 & -- & -- & -- & --\\
	DR& 0.02000 & 2.20 & 4 & 0 & -- & -- & -- & --\\
	
	\midrule[0.3pt]\bottomrule[1pt]
\end{tabular}
\label{tab:fuel24}
\end{table}

\begin{table}[htbp]
\centering\small
	\caption{Emission coefficients of thermal units}
\begin{tabular}[center]{c c c c c c p{.38cm} p{.38cm} p{.38cm} c c c}
	\toprule[1pt]\midrule[0.3pt]
	Unit & $\alpha_n$ & $\beta_n$ & $\gamma_n$\\ \midrule
	
	1 & 0.00312 & -0.24444 & 10.33908\\
	2 & 0.00312 & -0.24444 & 10.33908\\
	3 & 0.00509 & -0.40695 & 30.03910\\
	4 & 0.00509 & -0.40695 & 30.03910\\
	5 & 0.00344 & -0.38132 & 32.00006\\
	
	\midrule[0.3pt]\bottomrule[1pt]
\end{tabular}
\label{tab:emission24}
\end{table}

\section{Conclusion\label{sec:Con}}
In this paper, the UCD problem is studied in microgrid environment while taking account of thermal units and DGs and considering DR and CET.
Firstly, the UCD problem is reformulated as the optimal switching problem of a hybrid system.
The reformulation enables employing recent advances in areas such as RL and DP and allows easy handling of time-dependent costs.
Theoretical results show the equivalence between the original UCD problem the reformulated problem.
Then, the CLHO algorithm is developed by using neural networks to approximate the optimal cost-to-go and employing convex optimization algorithm to solve the lower problem.
The CLHO algorithm is capable of producing the optimal schedule in a closed-loop style, leading to the closed-loop solution for different initial conditions and external disturbances.
Theoretical analysis and simulation results are provided to show the effectiveness of the proposed algorithm.
It is also shown that the low emission quotas and high CET price help reduce the total carbon emissions.

\bibliographystyle{IEEEtran}
% Generated by IEEEtran.bst, version: 1.14 (2015/08/26)

\end{document}